%% file: main.tex
\documentclass{article}
\usepackage{iclr2025_conference,times}
\usepackage{amsthm,amssymb,amsmath}

\usepackage[utf8]{inputenc} 
\usepackage[T1]{fontenc}    
\usepackage{hyperref}       
\usepackage{url}            
\usepackage{booktabs}       
\usepackage{amsmath}       
\usepackage{nicefrac}       
\usepackage{microtype}      
\usepackage{wrapfig}

\usepackage{lipsum}         
\usepackage{graphicx}
\usepackage{natbib}
\usepackage{doi}

\usepackage{mathtools}
\usepackage{algorithm}
\usepackage{algpseudocode}
\usepackage{color,xcolor}
\usepackage{cleveref}       
\usepackage{etoolbox}
\usepackage{threeparttable} 

\input{cmd.tex}

\title{Decentralized Optimization with Coupled Constraints}

\author{
    \hspace{-0.08cm}Demyan Yarmoshik \\
	MIPT; Research Center for Artificial Intelligence, Innopolis University, Innopolis, Russia \\
	\texttt{yarmoshik.dv@phystech.edu}
	\And
	Alexander Rogozin \\
	  MIPT; Skoltech \\
	\texttt{aleksander.rogozin@phystech.edu}
	\AND
	Nikita Kiselev \\
	  Research Center for Artificial Intelligence, Innopolis University, Innopolis, Russia; MIPT \\
	\texttt{kiselev.ns@phystech.edu}
	\And
	  Daniil Dorin \\
	  MIPT \\
	\texttt{dorin.dd@phystech.edu}
	\AND
	  Alexander Gasnikov \\
	  Research Center for Artificial Intelligence, Innopolis University, Innopolis, Russia; MIPT; Skoltech \\
	\texttt{gasnikov@yandex.ru}
	\And
	Dmitry Kovalev \\
	Yandex Research \\
	\texttt{dakovalev1@gmail.com}
}

\iclrfinalcopy
\begin{document}
\maketitle

\begin{abstract}
We consider the decentralized minimization of a separable objective $\sum_{i=1}^{n} f_i(x_i)$, where the variables are coupled through an affine constraint $\sum_{i=1}^n\left(\mathbf{A}_i x_i - b_i\right) = 0$.
We assume that the functions $f_i$, matrices $\mathbf{A}_i$, and vectors $b_i$ are stored locally by the nodes of a computational network, and that the functions $f_i$ are smooth and strongly convex. 

This problem has significant applications in resource allocation and systems control and can also arise in distributed machine learning.
We propose lower complexity bounds for decentralized optimization problems with coupled constraints and a first-order algorithm achieving the lower bounds. To the best of our knowledge, our method is also the first linearly convergent first-order decentralized algorithm for problems with general affine coupled constraints.
\end{abstract}

\section{Introduction}\label{sec:intro}
We consider the decentralized optimization problem with coupled constraints 
\begin{equation}\label{eq:main}
  \min_{x_1 \in \R^{d_1},\ldots,x_n \in \R^{d_n}} \sum_{i=1}^{n} f_i(x_i) \quad \text{s.t.} \quad \sum_{i=1}^{n} (\mA_i x_i - b_i) = 0,
\end{equation}
where for $i \in \rng{1}{n}$ functions $f_i(x_i) \colon \R^{d_i} \to \R$ are continuously differentiable,
$\mA_i \in \R^{m \times d_i}$ and $b_i \in \R^m$  are constraint matrices and vectors respectively.

We are interested in solving problem~\eqref{eq:main} in a decentralized distributed setting. That is, we assume the existence of a communication network $\cG = (\cV, \cE)$, where $\cV = \rng{1}{n}$ is the set of compute nodes, and $\cE \subset \cV \times \cV$ is the set of communication links in the network. Each compute node $i \in \cV$ locally stores the objective function $f_i(x_i)$, the constraint matrix $\mA_i$ and the vector $b_i$. 
Compute node $i \in \cV$ can send information (e.g.,~vectors, scalars, etc.) to compute node $j \in \cV$ if and only if there is an edge $(i, j) \in \cE$ in the communication network.

Coupled constraints arise in various application scenarios, where sharing resources or information takes place. 
Often, due to the distributed nature of such problems, decentralization is desired for communication and/or privacy related reasons.
Let us briefly describe several practical cases of optimization problems with coupled constraints.

$\bullet$ \textbf{Optimal exchange.} 
Also known as the resource allocation problem~\cite{boyd2011distributed,nedic2018improved}, it writes as
\begin{align*}
  \min_{x_1,\ldots,x_n \in \sX} \sum_{i=1}^n f_i(x_i) \qst \sum_{i=1}^n x_i = b,
\end{align*}
where $x_i \in \sX$ represents the quantities of commodities exchanged among the agents of the system,
and $b \in \sX$ represents the shared budget or demand for each commodity.
This problem is essential in economics~\cite{arrow1954existence}, and systems control~\cite{dominguez2012decentralized}.

$\bullet$  \textbf{Problems on graphs.} 
In various applications, distributed systems are formed on the basis of physical networks. 
This is the case for electrical microgrids, telecommunication networks and drone swarms. 
Distributed optimization on graphs applies to such systems and encompasses, to name a few, 
optimal power flow~\cite{wang2016fully} 
and power system state estimation~\cite{zhang2024geometric} 
problems.

As an example, consider an electric power network.
Let $x_i \in \R^2$ denote the voltage phase angle and the magnitude at $i$-th electric node,
and let $s$ be the vector of (active and reactive) power flows for each pair of adjacent electric nodes.
Highly accurate linearization approaches~\cite{yang2016state, van2014dc} allow to formulate the necessary 
relation between voltages and power flows as a linear system of equations $\sum_{i=1}^n \mA_i x_i = s$.
An important property of the matrices $\mA_i$ is that their compatibility with the physical network (but not necessary with the communication network).
This means that for each row of the matrix $\cbr{\mA_1, \ldots, \mA_n}$, there is a node $k$ such that $\mA_i$ can have nonzero elements in this row only if nodes $i$ and $k$ are connected in the physical network, or $k = i$.


$\bullet$ \textbf{Consensus optimization.}
Related to the previous example is the consensus optimization~\cite{boyd2011distributed}
\begin{align*}
  \min_{x_1,\ldots,x_n \in \sX} \sum_{i=1}^n f_i(x_i) \qst x_1 = x_2 = \ldots = x_n.
\end{align*}
It is widely used in horizontal federated learning~\cite{kairouz2021advances}, as well as in the more general context of decentralized optimization of finite-sum objectives~\cite{gorbunov2022recent, scaman2017optimal}. 

To handle the consensus constraint, decentralized algorithms either reformulate it as $\sum_{i=1}^n \mW_i x_i = 0$, where $\mW_i$ is the $i$-th vertical block of a gossip matrix (an example of which is the communication graph's Laplacian), or utilize the closely related mixing matrix approach~\cite{gorbunov2022recent}.
Mixing and gossip matrices are used because they are communication-friendly: calculating the sum $\sum_{i=1}^n \mW_i x_i$ only requires each compute node to communicate once with each of its adjacent nodes.
Clearly, consensus optimization with gossip matrix reformulation can be reduced to~\eqref{eq:main} by setting $\mA_i = \mW_{i}$.
However, the principal difference between this example and~\eqref{eq:main}, is that~\eqref{eq:main} does not assume $\mA_i$ to be communication-friendly.
We discuss the complexity of the reduction from consensus optimization to optimization with coupled constraints in \Cref{sec:consensus_vs_coupled}.

$\bullet$ \textbf{Vertical federated learning (VFL).} 
In the case of VFL, the data is partitioned by features, differing from the usual (horizontal) federated learning, where the data is partitioned by samples \cite{yang2019federated, boyd2011distributed}.
Let $\mF$ be the matrix of features, split vertically between compute nodes into submatrices $\mF_i$, so that each node possesses its own subset of features for all data samples. 
Let $l \in \sY$ denote the vector of labels, and let $x_i \in \R^{d_i}$ be the vector of model parameters owned by the $i$-th node.
VFL problem formulates as
\begin{align}\label{eq:vfl-1}
  \min_{\substack{z \in \sY \\ x_1 \in \R^{d_1},\ldots,x_n \in \R^{d_n}}} \ell(z, l) + \sum_{i=1}^n r_i(x_i) \qst \sum_{i=1}^n \mF_i x_i = z,
\end{align}
where $\ell$ is a loss function, and $r_i$ are regularizers.
The constraints in~\eqref{eq:vfl-1} are coupled constraints, and the objective is separable; therefore, it is a special case of~\eqref{eq:main}.
We return to the VFL example in Section~\ref{sec:experiments}.

\noindent\textbf{Paper organization}. In Section~\ref{sec:related_work} we present a literature review. Subsequently, in Section~\ref{sec:mathematical_setting} we introduce the assumptions and problem parameters. Section~\ref{sec:derivation_of_algorithm} describes the key ideas of algorithm development and Section~\ref{sec:main_results} presents the convergence rate of the method and the lower complexity bounds. Finally, in Section~\ref{sec:experiments}, we provide numerical simulations.

\section{Related work and our contribution}\label{sec:related_work}
Decentralized optimization algorithms were initially proposed for consensus optimization~\cite{nedic2009distributed}, based on earlier research in distributed optimization~\cite{tsitsiklis1984problems, bertsekas1989parallel} and algorithms for decentralized averaging (\textit{consensus} or \textit{gossip} algorithms) \cite{boyd2006randomized, olshevsky2009convergence}, which assumed the existence of a communication network, as does the present paper.
The optimal complexity for consensus optimization was first achieved with a dual accelerated gradient descent in~\cite{scaman2017optimal}, where the method required computing gradients of Fenchel conjugates of $f_i(x)$. The corresponding complexity lower bounds were also established in the same paper. This result was later generalized to primal algorithms (which use gradients of the functions $f_i(x)$ themselves) \cite{kovalev2020optimal}, time-varying communication graphs~\cite{li2021accelerated,kovalev2021lower} and methods that use stochastic gradients~\cite{dvinskikh2019decentralized}.
Today there also exist algorithms with communication compression~\cite{beznosikov2023biased}, asynchronous algorithms~\cite{koloskova2024optimization}, algorithms for saddle-point formulations~\cite{rogozin2021decentralized} and gradient-free oracles~\cite{beznosikov2020derivative}, making decentralized consensus optimization a quite well-developed field~\cite{nedic2020distributed, gorbunov2022recent}, benefiting systems control~\cite{ram2009distributed} and machine learning~\cite{lian2017can}.

Beginning with the addition of local constraints to consensus optimization~\cite{nedic2010constrained, zhu2011distributed}, constrained decentralized optimization has been established as a research direction.
A zoo of distributed problems with constraints was investigated in~\cite{necoara2011parallel,necoara2014distributed, necoara2015linear}.

Primarily motivated by the demand from the power systems community, various decentralized algorithms for coupled constraints have been proposed. 
Generally designed for versatile engineering applications, many of these algorithms assume restricted function domains~\cite{wang2022distributed,liang2019distributed,nedic2018improved,gong2023decentralized,zhang2021distributed, wu2022distributed}, nonlinear inequality constraints~\cite{liang2019distributed,gong2023decentralized,wu2022distributed}, time-varying graphs~\cite{zhang2021distributed, nedic2018improved} or utilize specific problem structure~\cite{wang2022distributed}.


\begin{wrapfigure}[14]{r}{9cm}
\vspace{-0.93cm}
\begin{minipage}{9cm}
\begin{table}[H]
\centering
\begin{threeparttable}
\caption{Comparison of algorithms for decentralized optimization with coupled constraints}
  \label{tab:main}
  \begin{tabular}{|l|l|l|}
    \hline
    Reference & Oracle & Rate \\ \hline
    \cite{doan2017distributed}~\tnote{$\dagger$} & First-order & Linear\\ \hline
    \cite{falsone2020tracking} & Prox &  Sub-linear\\ \hline
    \cite{wu2022distributed} & Prox &  Sub-linear\\ \hline
    \cite{chang2016proximal}~ & Prox &  Sub-linear\\ \hline
    \cite{li2018accelerated}~\tnote{$\dagger$} & Prox & Linear \\ \hline
    \cite{gong2023decentralized} & Inexact prox & Linear \\ \hline
    \cite{nedic2018improved}~\tnote{$\dagger$} & First-order & Accelerated \\ \hline
    This work & First-order& Optimal \\ \hline
  \end{tabular}
  \begin{tablenotes}
    \footnotesize
    \item[$\dagger$] Applicable only for resource allocation problem
  \end{tablenotes}
\end{threeparttable}
\end{table}
\end{minipage}
\end{wrapfigure}

Works of \cite{doan2017distributed,li2018accelerated,nedic2018improved} focus on the resource allocation problem. 
For undirected time-varying graphs~\cite{doan2017distributed} proposes a first-order algorithm with $O(B\kappa_fn^2\ln\frac1\eps)$ communication and derivative computation complexity bound, where $B$ is the time required for the time-varying graph to reach connectivity.
\cite{li2018accelerated} applies a combination of gradient tracking and push-sum approaches from~\cite{nedic2017achieving} to obtain linear convergence on directed time-varying graphs in the restricted domain case, \textit{i.e.}, $x_i \in \Omega_i$, where $\Omega_i$ is a nonempty closed convex set. 
\cite{nedic2018improved} achieves accelerated linear convergence via a proximal point method in the restricted domain case.
When $\Omega_i = \R^d$, they also show that Nesterov's accelerated gradient descent can be applied to achieve optimal $O(\sqrt{\kappa_\mW}\sqrt{\kappa_f}\ln\frac1\eps)$ communication complexity.
In~\cite{gong2023decentralized} an inexact proximal-point method is proposed to solve problems with coupled affine equality and convex inequality constraints.
Linear convergence is proved when the inequalities are absent, and $\Omega_i$ are convex polyhedrons. The papers \cite{wu2022distributed}, \cite{chang2016proximal}, \cite{falsone2020tracking} present algorithms with sub-linear convergence.

As summarized in Table~\ref{tab:main}, no accelerated linearly convergent algorithms for general affine-equality coupled constraints were present in the literature prior to our work.
Also, most of the algorithms require proximal oracle, which allows to handle more general problem formulations, but has higher computational burden than the first-order oracle.
We propose a new first-order decentralized algorithm with optimal (accelerated) linear convergence rate. 
We prove its optimality by providing lower bounds for the number of objective's gradient computations, matrix multiplications and decentralized communications, which match complexity bounds for our algorithm.

\section{Mathematical setting and assumptions}\label{sec:mathematical_setting}

Let us begin by introducing the notation. 
The largest and smallest nonzero eigenvalues (or singular values) of a matrix $\mC$ are denoted by $\lmax(\mC)$~ (or $\smax(\mC)$) and $\lminp(\mC)$~ (or $\sminp(\mC)$), respectively. For vectors $x_i \in \R^{d_i}$ we introduce a column-stacked vector $x = \col(x_1, \ldots, x_m) = (x_1^\top \ldots x_m^\top)^\top\in\R^{d}$. We denote the identity matrix by $\mI_m \in \R^{m \times m}$. The symbol $\otimes$ denotes the Kronecker product of matrices. 
By $\cL_m$ we denote the so-called consensus space, which is given as $\cL_m = \{(y_1,\ldots,y_n) \in (\R^m)^n : y_1 , \ldots , y_n \in \R^m\; \text{ and }\; y_1 = \cdots = y_n\}$, and $\cL_m^\perp$ denotes the orthogonal complement to $\cL_m$, which is given as
\begin{equation}\label{eq:Lperp}
	\cL_m^\perp = \{(y_1,\ldots,y_n) \in (\R^m)^n : y_1 , \ldots , y_n \in \R^m\;\text{ and }\; y_1 + \cdots + y_n = 0\}.
\end{equation}

\begin{assumption}\label{ass:f}
  Continuously differentiable functions $f_i(x) \colon \R^{d_i} \to \R , \ i \in \rng{1}{n} $ are $L_f$-smooth and $\mu_f$-strongly convex, where $L_f \geq \mu_f > 0$.
  That is, for all $x_1,x_2 \in \R^{d_i}$ and $i \in \rng{1}{n}$, the following inequalities hold:
	\begin{equation*}
		\frac{\mu_f}{2}\sqn{x_2 - x_1}
		\leq
		f_i(x_2) - f_i(x_1) - \<\nabla f_i(x_1), x_2 - x_1>
		\leq
		\frac{L_f}{2}\sqn{x_2 - x_1}.
	\end{equation*}
	By $\kappa_{f}$ we denote the condition number $\kappa_{f} = {L_{f}}/{\mu_{f}}$.
\end{assumption}

\begin{assumption}\label{ass:A}
    There exists $x^* = (x^*_1,\ldots,x^*_n), x^*_i \in \R^{d_i}$ such that $\sum_{i=1}^{n} (\mA_i x^*_i - b_i) = 0$.
	There exist constants $L_{\mA}\geq \mu_{\mA} > 0$, such that the constraint matrices  $\mA_1,\ldots,\mA_n$ satisfy the following inequalities:
	\begin{equation}\label{eq:def_mu_a_L_a}
		\smax^2(\mA) = \max_{i \in \rng{1}{n}} \smax^2(\mA_i) \leq L_{\mA},
    \qquad
		\mu_{\mA} \leq \lminp\left( \mS \right),
	\end{equation}
	where the matrix $\mS \in \R^{m\times m}$ is defined as $\mS = \frac{1}{n}\sum_{i=1}^{n} \mA_i \mA_i^\top$. 
    We also define the condition number of 
the block-diagonal matrix $\mA = \diag{\mA_1, \ldots,\mA_n}\in \R^{mn \times d}$
 as $\kappa_{\mA} = L_\mA / \mu_\mA$.
\end{assumption}

For any matrix $\mM$ other than $\mA$ we denote by $L_\mM$ and $\mu_\mM$ some upper and lower bound on its maximal and minimal positive squared singular values respectively: 
\begin{equation}\label{eq:mat_Lmu}
  \lmax(\mM^\top\mM) = \smax^2(\mM) \leq L_{\mM},
    \qquad
	\mu_{\mM} \leq \sminp^2(\mM) = \lminp(\mM^\top \mM).
\end{equation}

We also assume the existence of a so-called gossip matrix $W \in \R^{n\times n}$ associated with the communication network $\cG$, which satisfies the following assumption.
\begin{assumption}\label{ass:W}
    The gossip matrix $W$ is a $n\times n$ symmetric positive semidefinite matrix such that:\\
    \hspace*{1em} 1. $W_{ij}\neq 0$ if and only if $(i, j)\in\cE$ or $i = j$.\\
    \hspace*{1em} 2. $Wy = 0$ if and only if $y\in\cL_1$, i.e. $y_1 = \ldots = y_n$.\\
    \hspace*{1em} 3. There exist constants $L_\mW\geq \mu_\mW > 0$ such that $\mu_\mW\leq \lminp^2(W)$ and $\lmax^2(W)\leq L_\mW$.
\end{assumption}
We will use a dimension-lifted analogue of the gossip matrix defined as $\mW = W \otimes \mI_m$.
From the properties of the Kronecker product of matrices it follows that $\lminp^2(\mW) = \lminp^2(W)$ and $\lmax^2(\mW) = \lmax^2(W)$. By $\kappa_{\mW}$ we denote the condition number
\begin{equation}\label{eq:kappa_mW}
   \kappa_{\mW} =  \sqrt\frac{L_{\mW}}{\mu_{\mW}} \geq \frac{\lmax(\mW)}{\lminp(\mW)}.
\end{equation}
Moreover, the kernel and range spaces of $W$ and $\mW$ are given by
\begin{align}
	\label{eq:range_W}
  \ker W = \cL_1, \ \range W = \cL_1^\perp, \quad \ker \mW = \cL_m, \ \range \mW = \cL_m^\perp.
\end{align}

\section{Derivation of the algorithm}\label{sec:derivation_of_algorithm}
\subsection{Strongly convex communication-friendly reformulation}
Let $\pW$ be any positive semidefinite matrix such that 
\begin{equation}\label{eq:range_pW}
  \range \pW = (\ker \pW)^\perp =  \cL_m^\perp,
\end{equation}
and multiplication of a vector $y = (y_1, \ldots, y_n) \in (\sY)^n$ by $\pW$ can be performed efficiently in the decentralized manner if its $i$-th block component $y_i$ is stored at $i$-th node of the computation network.
Similarly to~\cref{eq:kappa_mW}, we define
\begin{equation}\label{eq:kappa_pW}
    \kappa_{\pW} = \sqrt\frac{L_{\pW}}{\mu_{\pW}} \geq \frac{\lmax(\pW)}{\lminp(\pW)}.
\end{equation}
Due to the definition of $\mW$ and \cref{eq:range_W}, the simplest choice for $\pW$ might be to set $\pW = \mW$. 
Later we will specify another way to choose $\pW$ for optimal algorithmic performance.

Problem~\eqref{eq:main} can be reformulated as follows:
\begin{equation}\label{eq:lifted}
	\min_{x\in \sX, y \in (\sY)^n} G(x,y) \quad\text{s.t.}\quad \mA x + \gamma \pW y = \mb,
\end{equation}
where the function $G(x,y) \colon \sX \times (\sY)^n \to \R$ is defined as
\begin{equation}\label{eq:G}
	G(x,y) = F(x) + \frac{r}{2}\sqn{\mA x + \gamma \pW y - \mb},
\end{equation}
the function $F(x) \colon \sX \to \R$ is defined as
$F(x) = \sum_{i=1}^{n} f_i(x_i)$, where $x=(x_1,\ldots,x_n), x_i \in \R^{d_i}$,
the matrix $\mA \in \R^{mn \times d}$ is the block-diagonal matrix $\mA = \diag{\mA_1, \ldots,\mA_n}$,
the vector $\mb$ is the column-stacked vector $\mb = \col\cbraces{b_1,\ldots,b_n} \in \R^{mn}$,
and $r, \gamma > 0$ are scalar constants that will be determined later.

From the definitions of $\mA$, $\mb$ and $\cL_m^\bot$ (\cref{eq:Lperp}) it is clear that $\sum_{i=1}^n (\mA_i x_i - b_i) = 0$ if and only if $\mA x - \mb \in \cL_{m}^\perp$. 
Since $\range \pW = \cL^\perp_m$, the constraint in problem~\eqref{eq:lifted} is equivalent to the coupled constraint in~\eqref{eq:main}.
For all $x, y$ satisfying the constraint, the augmented objective function $G(x, y)$ is equal to the original objective function $F(x)$. 
Therefore, problem \eqref{eq:lifted} is equivalent to problem \eqref{eq:main}.


The following \Cref{lem:G} shows that the function $G(x,y)$ is strongly convex and smooth.

\begin{lemma}\label{lem:G}
	Let $r$ and $\gamma$ be defined as follows:
	\begin{equation}\label{eq:r_gamma}
		r = \frac{\mu_f}{2L_{\mA}}, \quad \gamma^2 = \frac{\mu_\mA + L_\mA}{\mu_\pW}.
	\end{equation}

  Then, the strong convexity and smoothness constants of $G(x, y)$ on $\sX \times \cL_m^\perp $ are given by 
	\begin{equation}\label{eq:G_mu_L}
    \mu_G = \mu_f \min\braces{\frac12, \frac{\mu_\mA + L_\mA}{4 L_\mA}}, \quad L_G = \max\braces{L_f + \mu_f,  \mu_f\frac{\mu_\mA + L_\mA}{L_\mA}\frac{L_\pW}{\mu_\pW}}.
	\end{equation}
\end{lemma}

Let the matrix $\mB \in \R^{mn \times (d+mn)}$ be defined as $\mB = \left[\begin{matrix} \mA & \gamma \pW \end{matrix}\right]$.
The following \Cref{lem:B} connects the spectral properties of $\mB$, $\mA$ and $\pW$.
\begin{lemma}\label{lem:B}
	The following bounds on the singular values of $\mB$ hold:
  \begin{equation}\label{eq:B_mu_L}
     \sminp^2(\mB) \geq \mu_\mB = \frac{\mu_\mA}{2} , \quad  \smax^2(\mB) \leq L_\mB = L_\mA + (L_\mA + \mu_\mA)\frac{L_\pW}{\mu_\pW},
  \end{equation}
  and
  \begin{equation}\label{eq:kappa_B}
    \frac{\smax^2(\mB)}{\sminp^2(\mB)} \leq \frac{L_\mB}{\mu_\mB} = \kappa_\mB = 2\cbraces{\kappa_\mA + \frac{L_\pW}{\mu_\pW}(1+\kappa_\mA)}.
  \end{equation}
\end{lemma}

Proofs of \Cref{lem:G} and \Cref{lem:B} are provided in Appendix~\ref{app:derivation_of_algorithm}.

\subsection{Chebyshev acceleration}

Chebyshev acceleration allows us to decouple the number of computations of the objective's gradient $\nabla F(x)$ from the properties of the communication network and the constraint matrix~--- specifically, from the condition numbers $\kappa_\mW$ and $\kappa_\mA$. The Chebyshev trick enables to replace the matrix with a matrix polynomial with a better condition number.

Consider some affine relation $\mM  u = \bd$ and let $\cP_\mM$ be a polynomial such that $\cP_\mM(\lambda) = 0 \Leftrightarrow \lambda = 0$ for any eigenvalue $\lambda$ of $\mM^\top\mM$. Note that here we interchangeably use $\cP$ as a polynomial of a matrix and a polynomial of a scalar. We denote any feasible point for the constraint $\mM u = \bd$ as $u_0$.
Then,
\begin{align*}
	\mM u = \bd 
	&\Leftrightarrow
	\mM (u - u_0) = 0
	\alra{is due to $\ker \mM^\top\mM = \ker \mM$}
	\mM^\top \mM (u - u_0) = 0
	\\& \alra{is due to $\ker \cP_\mM(\mM^\top \mM) = \ker \mM^\top \mM$ by the assumption about $\cP_\mM(\lambda)$}
	\cP_\mM(\mM^\top \mM) (u - u_0) = 0 
	\alra{is due to $\ker \mM^\top\mM = \ker \mM$}
	\sqrt{\cP_\mM(\mM^\top \mM)} (u - u_0) = 0 
\end{align*}
where \annotate. 

%

Following~\cite{salim2022optimal} and~\cite{scaman2017optimal}, we use the translated and scaled Chebyshev polynomials, because they are the best at compressing the spectrum \cite{auzinger2011iterative}.
 
\begin{lemma}[\cite{salim2022optimal}, Section 6.3.2]\label{lem:poly}
  Consider a matrix $\mM$.
  Let $\ell = \left\lceil\sqrt{\frac{L_\mM}{\mu_\mM}}\right\rceil \geq \left\lceil\sqrt{\frac{\lmax(\mM^\top\mM)}{\lminp(\mM^\top\mM)}}\right\rceil$.
  Define $\cP_{\mM}(t) = 1 - \frac{T_\ell\cbraces{(L_\mM + \mu_\mM - 2t) / (L_\mM - \mu_\mM)}}{T_\ell\cbraces{(L_\mM + \mu_\mM) / (L_\mM - \mu_\mM)}}$,
  where $T_\ell$ is the Chebyshev polynomial of the first kind of degree $n$ defined by $T_\ell(t) = \frac12 \cbraces{\cbr{t + \sqrt{t^2 - 1}}^\ell + \cbr{t - \sqrt{t^2 - 1}}^\ell }$.
  Then, $\cP_{\mM}(0) = 0$, and 
  \begin{align}
    \lmax\cbraces{\cP_{\mM}(\mM^\top\mM)} &\leq \max_{t \in [\mu_\mM, L_\mM]} \cP_{\mM}(t) \leq \frac{19}{15},
   \\
    \lminp\cbraces{\cP_{\mM}(\mM^\top\mM)} &\geq \min_{t \in [\mu_\mM, L_\mM]} \cP_{\mM}(t) \geq \frac{11}{15}.
  \end{align}
\end{lemma}

Results of this section are summarized in the following \Cref{lem:prob_cheb}.
\begin{lemma}\label{lem:prob_cheb}
  Define
  \begin{equation}\label{eq:Wcheb}
    \pW = \cP_{\sqrt\mW}(\mW)
  \end{equation}
  and
  \begin{equation}\label{eq:K}
    \mK = \sqrt{\cP_{\mB}(\mB^\top\mB)}.
  \end{equation}
Let $G(u) = G(x, y)$, 
$\sU = \sX \times \cL_m^\perp$ and
$\pb = \sqrt{\cP_{\mB}(\mB^\top \mB)}u_0$.
Then, problem
\begin{equation}\label{eq:compact}
	\min_{u\in \sU} G(u) \quad\text{s.t.}\quad \mK u = \pb
\end{equation}
is an equivalent preconditioned reformulation of problem~\eqref{eq:lifted}, and, in turn, of problem~\eqref{eq:main}. 
\end{lemma}

\subsection{Base algorithm}
\begin{wrapfigure}[12]{r}{9cm}
\vspace{-0.45cm}
\begin{minipage}{9cm}
  \begin{algorithm}[H]
      \caption{APAPC}
  	\label{alg:apapc}
  	\begin{algorithmic}[1]
  		\State {\bf Parameters:}  
      $u^0 \in \cU$
  		$\eta,\theta,\alpha>0$, $\tau \in (0,1)$
  		\State Set $u_f^0 = u^0$, $z^0 = 0 \in \cU$
  		\For{$k=0,1,2,\ldots$}{}
  		\State $u_g^k \eqdef \tau u^k +  (1-\tau)u_f^k$\label{alg:apapc:line:x:1}
  		\State $u^{k+\frac{1}{2}} \eqdef (1+\eta\alpha)^{-1}(u^k - \eta (\nabla G(u_g^k) - \alpha  u_g^k + z^k))$\label{alg:apapc:line:x:2}
  		\State $z^{k+1} \eqdef z^k  + \theta \mK^\top (\mK u^{k+\frac{1}{2}} - \pb)$ \label{alg:apapc:line:z}
  		\State $u^{k+1} \eqdef (1+\eta\alpha)^{-1}(u^k - \eta (\nabla G(u_g^k) - \alpha  u_g^k + z^{k+1}))$\label{alg:apapc:line:x:3}
  		\State $u_f^{k+1} \eqdef u_g^k + \tfrac{2\tau}{2-\tau}(u^{k+1} - u^k)$\label{alg:apapc:line:x:4}
  		\EndFor
  	\end{algorithmic}
  \end{algorithm}
\end{minipage}
\end{wrapfigure}

Our base algorithm, \Cref{alg:apapc}, is the Proximal Alternating Predictor-Corrector (PAPC) with Nesterov's acceleration, called Accelerated PAPC (APAPC).
It was proposed in~\cite{salim2022optimal} to obtain an optimal algorithm for optimization problems formulated as~\eqref{eq:compact}.
See~\cite{kovalev2020optimal, salim2022dualize} for the review of related algorithms and history of their development.

APAPC algorithm formulates as \Cref{alg:apapc}, and its convergence properties are given in \Cref{prop:apapc}.

\begin{proposition}[\cite{salim2022optimal}, Proposition 1]\label{prop:apapc}
  Assume that the matrix $\mK$ in \eqref{eq:compact} satisfies $\mu_{\mK} > 0$ and $\mb' \in \range \mK$, and denote $\kappa_\mK = \frac{L_\mK}{\mu_\mK}$.
  Also assume that the function $G$ is $L_G$-smooth and $\mu_G$-strongly convex.
  Set the parameter values of \Cref{alg:apapc} as $\tau = \min\braces{1, \frac12\sqrt\frac{\kappa_\mK}{\kappa_G}}$, $\eta= \frac1{4\tau L_G}$, $\theta = \frac1{\eta L_\mK}$ and $\alpha = \mu_G$.
  Denote by $u^*$ the solution of problem~\eqref{eq:compact} and by $z^*$ the solution of its dual problem satisfying $z^* \in \range \mK$.
  Then the iterates $u^k, z^k$ of \Cref{alg:apapc} satisfy
    \begin{align}
		&\frac{1}{\eta}\sqN{u^k - u^{\star}} + \frac{\eta\alpha}{\theta(1+\eta\alpha)}\sqN{(\mK^\top)^\dagger z^k - z^{\star}} \\
		&+
			\frac{2(1-\tau)}{\tau}\bg_G(u_f^k, u^{\star})
			\leq
			\cbraces{\!1 + \frac{1}{4}\min\braces{\frac{1}{\sqrt{\kappa_G\kappa_\mK}},\frac{1}{\kappa_\mK}}\!}^{-k}
			C,\notag
		\end{align}
		where
		$C \eqdef \frac{1}{\eta}\sqN{u^0 - u^{\star}} + \frac{1}{\theta}\sqn{z^0 - z^{\star}} +
		\frac{2(1-\tau)}{\tau}\bg_G(u_f^0, u^{\star})$,
    and $\bg_G$ denotes the Bregman divergence of $G$, defined by $D_G(u', u) = G(u') - G(u) - \<\nabla G(u), u' - u>$.
\end{proposition}

\section{Main results}\label{sec:main_results}
\subsection{Algorithm}\label{subsec:algorithm}
\begin{wrapfigure}[15]{r}{7.8cm}
\begin{minipage}{7.8cm}
  \vspace{-0.87cm}
  \begin{algorithm}[H]
      \caption{Main algorithm}
  	\label{alg:main}
  	\begin{algorithmic}[1]
  		\State {\bf Parameters:}  
      $x^0 \in \sX$
  		$\eta,\theta,\alpha>0$, $\tau \in (0,1)$
  
  		\State Set $y^0 \eqdef 0 \in (\sY)^n$, $u^0 \eqdef (x^0, y^0)$,
      \\  \quad \ \ $u_f^0 \eqdef u^0$, $z^0 \eqdef 0 \in \sX \times (\sY)^n$
  		\For{$k=0,1,2,\ldots$}{}
  		\State $u_g^k \eqdef \tau u^k +  (1-\tau)u_f^k$
      \State $g^k \eqdef \gradG(u_g^k) - \alpha u_g^k$
  		\State $u^{k+\frac{1}{2}} \eqdef (1+\eta\alpha)^{-1}(u^k - \eta (g^k + z^k))$
      \State $z^{k+1} \eqdef z^k  + \theta \cdot \Kcheb(u^{k+\frac12})$
  		\State $u^{k+1} \eqdef (1+\eta\alpha)^{-1}(u^k - \eta (g^k + z^{k+1}))$
  		\State $u_f^{k+1} \eqdef u_g^k + \tfrac{2\tau}{2-\tau}(u^{k+1} - u^k)$
  		\EndFor
  	\end{algorithmic}
  \end{algorithm}
\end{minipage}
\end{wrapfigure}
As stated in \Cref{lem:prob_cheb}, problem~\eqref{eq:compact} is equivalent to problem~\eqref{eq:main}.
Due to~\Cref{lem:G}, its objective is strongly convex, allowing us to apply \Cref{alg:apapc} to it.
Using~\Cref{lem:poly}, we obtain that the condition numbers of $\pW$ and $\mK$ are bounded as $O(1)$,
but a single multiplication by $\pW$ and $\mK, \mK^\top$ translates to $O(\sqrt{\kappa_\mW})$ multiplications by $\mW$ and $O(\sqrt{\kappa_\mB})$ multiplications by $\mB, \mB^\top$ respectively.

We implement multiplications by $\pW$ and $\mK, \mK^\top$ through numerically stable Chebyshev iteration procedures given in \Cref{alg:mult_pW,alg:Kcheb}, which only use decentralized communications and multiplications by $\mA, \mA^\top$. 
\Cref{lem:G,lem:B,lem:poly} allow us to express the complexity of \Cref{alg:apapc} in terms of the parameters of the initial problem given in \Cref{ass:f,ass:A,ass:W}.
All this leads us to the following \Cref{thm:main}, a detailed proof of which is provided in Appendix~\ref{app:thm:main}, as well as the derivation of \Cref{alg:mult_pW,alg:Kcheb} and values of the parameters of \Cref{alg:main}.


\begin{figure*}[t!]
\begin{minipage}{.482\textwidth}
  \begin{algorithm}[H]
    \caption{$\mulpW(y)$: Multiplication by $\pW$}
  	\label{alg:mult_pW}
  	\begin{algorithmic}[1]
  		\State {\bf Parameters:} $y$ 
      \State $\rho\eqdef \cbraces{\sqrt{L_\mW}-\sqrt{\mu_\mW}}^2/16$
      \State $\nu\eqdef\cbraces{\sqrt{L_\mW}+\sqrt{\mu_\mW}}/2$
  		\State $\delta^0 \eqdef -\nu/2$, \ $n \eqdef \ceil{\sqrt{\kappa_\mW}}$
      \State $p^0 \eqdef -\mW y/\nu$, \ $y^1 \eqdef y + p^0$
  		\For{$i=1,\ldots,n-1$}{}
  		\State $\beta^{i-1}\eqdef\rho/\delta^{i-1}$
  		\State $\delta^i\eqdef-(\nu+\beta^{i-1})$
  		\State $p^i\eqdef\cbraces{\mW y^{i}+\beta^{i-1}p^{i-1}}/\delta^i$
  		\State $y^{i+1}\eqdef y^i+p^i$
  		\EndFor
  	\State {\bf Output:} $y - y^n$
  	\end{algorithmic}
  \end{algorithm}
  \vspace{-0.7cm}
  \begin{algorithm}[H]
    \caption{$\gradG(u)$: Computation of $\nabla G(u)$}
  	\label{alg:gradG}
  	\begin{algorithmic}[1]
  		\State {\bf Parameters:} $u = (x, y)$ 
      \State $z \eqdef r\cbraces{\mA x + \gamma \cdot\mulpW(y) - \mb}$
        \State {\bf Output:} $\pmat{\nabla F(x) + \mA^\top z \\ \gamma\cdot \mulpW(z)}$
  	\end{algorithmic}
  \end{algorithm}

\end{minipage}\ \ \ \
\begin{minipage}{.482\textwidth}
    \vspace{0.2cm}
  \begin{algorithm}[H]
    \caption{$\Kcheb(u)$: Computation of $\mK^\top(\mK u - \pb)$}
  	\label{alg:Kcheb}
  	\begin{algorithmic}[1]
  		\State {\bf Parameters:} $u = (x, y)$ 
      \State $\rho\eqdef \cbraces{{L_\mB}-{\mu_\mB}}^2/16$
      \State $\nu\eqdef\cbraces{{L_\mB}+{\mu_\mB}}/2$
  		\State $\delta^0 \eqdef -\nu/2$, \ $n \eqdef \ceil{\sqrt{\kappa_\mB}}$
      \State $q^0 \eqdef \mA x + \gamma \cdot\mulpW(y) -\mb$
      \State $p^0 \eqdef -\dfrac1\nu \pmat{\mA^\top  q^0 \\ \gamma \cdot \mulpW(q^0)}$
  		\State $u^1 \eqdef u + p^0$
  		\For{$i=1,\ldots,n-1$}{}
  		\State $\beta^{i-1}\eqdef\rho/\delta^{i-1}$
  		\State $\delta^i\eqdef-(\nu+\beta^{i-1})$
      \State $(x^i, y^i) = u^i$
      \State $q^i \eqdef \mA x^i + \gamma \cdot\mulpW(y^i) -\mb$
      \State $p^i\eqdef \dfrac1{\delta^i}\pmat{\mA^\top  q^i \\ \gamma \cdot \mulpW(q^i)}+\beta^{i-1}p^{i-1}/\delta^i$
  		\State $u^{i+1}\eqdef u^i+p^i$
  		\EndFor
  	\State {\bf Output:} $u - u^n$
  	\end{algorithmic}
  \end{algorithm}
\end{minipage}
\end{figure*}

\begin{theorem}\label{thm:main} Set the parameter values of \Cref{alg:main} as $\tau = \min\braces{1, \frac12\sqrt\frac{19}{44 \max\braces{1+\kappa_f, 6}}}$, $\eta= \frac1{4\tau \max\braces{L_f + \mu_f, 6\mu_f}}$, $\theta = \frac{15}{19\eta }$ and $\alpha = \frac{\mu_f}{4}$.
  Denote by $x^*$ the solution of problem~\eqref{eq:main}.
  Then, for every $\eps > 0$, \Cref{alg:main} finds $x^k$ for which $\sqn{x^k - x^*} \leq \eps$ using $O(\sqrt{\kappa_f}\log(1/\eps))$ objective's gradient computations, $O(\sqrt{\kappa_f}\sqrt{\kappa_\mA}\log(1/\eps))$ multiplications by $\mA$ and $\mA^\top$, and $O(\sqrt{\kappa_f}\sqrt{\kappa_\mA}\sqrt{\kappa_\mW}\log(1/\eps))$ communication rounds (multiplications by $\mW$).
\end{theorem}

\subsection{Lower bounds}\label{subsec:lower_bounds}

Let us formulate the lower complexity bounds for decentralized optimization with affine constraints. To do that, we formalize the class of the algorithms of interest. In the literature, approaches with continuous time~\cite{scaman2017optimal} and discrete time~\cite{kovalev2021lower} are used. We use the latter discrete time formalization. We assume that the method works in synchronized rounds of three types: local objective's gradient computations, local matrix multiplications and communications. At each time step, algorithm chooses one of the three step types.

Since the devices may have different dimensions $d_i$ of locally held vectors $x_i$, they cannot communicate these vectors directly. Instead, the nodes exchange quantities $\mA_i x_i\in\R^m$. For this reason, we introduce two types of memory $\cM_i(k)$ and $\cH_i(k)$ for node $i$ at step $k$. Set $\cM_i(k)$ stands for the local memory that the node does not share and $\cH_i(k)$ denotes the memory that the node exchanges with neighbors. The interaction between $\cM_i(k)$ and $\cH_i(k)$ is performed via multiplications by $\mA_i$ and $\mA_i^\top$.

Memory is initialized as $\cM_i(0) = \braces{0},~ \cH_i(0) = \braces{0}$. Below we describe how the sets $\cM_i(k), \cH_i(k)$ are updated.

1. Algorithm performs local gradient comutation round at step $k$. Gradient updates only operate in $\cM_i(k)$ and do not affect $\cH_i(k)$. For all $i\in\cV$ we have
\begin{align*}
	\cM_i(k+1) = \Span\braces{x, \nabla f_i(x), \nabla f_i^*(x):~ x\in M_i(k)},~ \cH_i(k+1) = \cH_i(k),
\end{align*}
where $f_i^*$ is the Fenchel conjugate of $f_i$.\\
2. Algorithm performs local matrix multiplication round at step $k$. Sets $\cH_i(k)$ and $\cM_i(k)$ make mutual updates via multiplication by $\mA_i$ and $\mA_i^\top$. For all $i\in\cV$ we have
\begin{align*}
	\cM_i(k+1) = \Span\braces{\mA_i^\top b_i,~ \mA_i^\top y:~ y\in \cH_i(k)},~ \cH_i(k+1) = \Span\braces{b_i, \mA_i x:~ x\in \cM_i(k)}.
\end{align*}
3. Algorithm performs a communication round at step $k$. The non-shared local memory $\cM_i(k)$ stays unchanged, while the shared memory $\cH_i(k+1)$ is updated via interaction with neighbors. For all $i\in\cV$ we have
\begin{align*}
	\cM_i(k+1) = \cM_i(k),~ \cH_i(k+1) = \Span\braces{\cH_j(k):~ (i, j)\in\cE}.
\end{align*}

Under given memory and computation model, we formulate the lower complexity bounds.
\begin{theorem}\label{thm:lower_bounds}
	For any $L_f > \mu_f > 0,~ \kappa_\mA, \kappa_\mW > 0$ there exist $L_f$-smooth $\mu_f$-strongly convex functions $\braces{f_i}_{i=1}^n$, matrices $\mA_i$ such that $\kappa_\mA = L_\mA / \mu_\mA$ (where $L_\mA, \mu_\mA$ are defined in \eqref{eq:def_mu_a_L_a}), and a communication graph $\cG$ with a corresponding gossip matrix $\mW$ such that $\kappa_\mW = \lambda_{\max}(\mW) / \lambda_{\min}^+(\mW)$, for which any first-order decentralized algorithm on problem \eqref{eq:main} to reach accuracy $\eps$ requires at least
	\begin{align*}
		N_f &= \Omega\cbraces{\sqrt\kappa_f \log\cbraces{\frac{1}{\eps}}}~ \text{gradient computations}, \\
		N_{\mA} &= \Omega\cbraces{\sqrt{\kappa_f}\sqrt{\kappa_\mA}\log\cbraces{\frac{1}{\eps}}}~ \text{multiplications by $\mA$ and $\mA^\top$}, \\
		N_{\mW} &= \Omega\cbraces{\sqrt{\kappa_f}\sqrt{\kappa_\mA}\sqrt{\kappa_\mW}\log\cbraces{\frac{1}{\eps}}}~ \text{communication rounds (multiplications by $\mW$)}.  \\
	\end{align*}
\end{theorem}
A proof of Theorem~\ref{thm:lower_bounds} is provided in Appendix~\ref{app:lower_bounds_proof}.

\section{Experiments}\label{sec:experiments}
The experiments were run on CPU Intel(R) Core(TM) i9-7980XE, with 62.5 GB RAM.

\begin{figure}[h!]
    \centering
    \includegraphics[width=\textwidth]{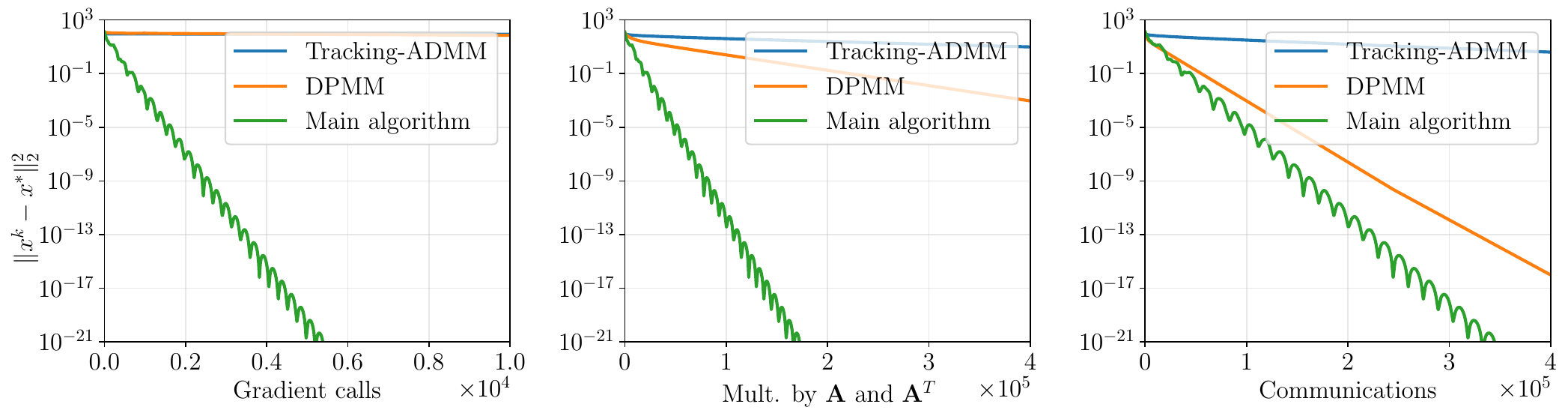}
    \vspace{-0.5cm}
    \caption{Synthetic, Erd\H{o}s–R\'enyi graph, $n = 20$, $d_i = 3$, $m = 10$}
    \label{fig:example-1}
\end{figure}

$\bullet$ \textbf{Synthetic linear regression.}
In this section we perform numerical experiments on a synthetic linear regression problem with $\ell_2$-regularization:
\begin{equation}\label{eq:example}
  \min_{x_1, \ldots, x_n \in \R^{d_i}} \sum_{i=1}^{n} \left( \frac{1}{2}\| C_i x_i - d_i \|_2^{2} + \frac{\theta}{2} \| x_i \|_2^2 \right) \quad \text{s.t.} \quad \sum_{i=1}^{n} (\mA_i x_i - b_i) = 0,
\end{equation}
where we randomly generate matrices $C_i \in \mathbb{R}^{d_i \times d_i}$, $\mathbf{A}_i \in \mathbb{R}^{m \times d_i}$ and vectors $d_i \in \mathbb{R}^{d_i}$, $b_i \in \mathbb{R}^{m}$ from the standard normal distribution. Local variables $x_i \in \mathbb{R}^{d_i}$ have the same dimension $d_i$, equal for all devices. Regularization parameter $\theta$ is $10^{-3}$. In the Fig.~\ref{fig:example-1} we demonstrate the performance of the our method on the problem, that has the following parameters: $\kappa_f = 3140$, $\kappa_{\mathbf{A}}=27$, $\kappa_{\mathbf{W}}=89$. There we use Erd\H{o}s–R\'enyi graph topology with $n = 20$ nodes. Local variables dimension is $d_i = 3$ and number of linear constraints is $m = 10$. 
We compare performance of \Cref{alg:main} with Tracking-ADMM algorithm \cite{falsone2020tracking} and DPMM algorithm \cite{gong2023decentralized}.
Note that Tracking-ADMM and DPMM are proximal algorithms that solve a subproblem at each iteration. The choice of objective function in our simulations (linear regression) makes the corresponding proximal operator effectively computable via Conjugate Gradient algorithm~\cite{nesterov2004introduction} that uses gradient computations. Therefore, we measure the computational complexity of these methods in the number of gradient computations, not the number of proximal operator computations.

$\bullet$ \textbf{VFL linear regression on real data.} 
Now we return to the problem, that we have announced in the introduction section. We apply VFL in the linear regression problem: $\ell$ is a typical mean squared loss function, that is $\ell(z, l) = \frac{1}{2} \| z - l \|_2^2$, and $r_i$ are $\ell_2$-regularizers, \textit{i.e.} $r_i(x_i) = \lambda \| x_i \|_2^2$. To adapt this from \eqref{eq:vfl-1} to \eqref{eq:main}, we redefine $x_1 := \binom{x_1}{z}$ and $x_2 := x_2, \ldots, x_n := x_n$. Thus, we can derive constraints matrices as in the \eqref{eq:main}:
\begin{align}\label{eq:vfl-matr}
    \mathbf{A}_1 = \begin{pmatrix} \mathbf{F}_1 & -\mathbf{I} \end{pmatrix}, \qquad \mathbf{A}_1 x_1 = \mathbf{F}_1 w_1 - z,\\
\mathbf{A}_i = \mathbf{F}_i, \quad i = 2, \ldots, n, \qquad \sum\limits_{i=1}^{n} \mathbf{A}_i x_i = \sum\limits_{i=1}^{n} \mathbf{F}_i w_i - z.
\end{align}

\begin{figure}[H]
    \centering
    \includegraphics[width=\textwidth]{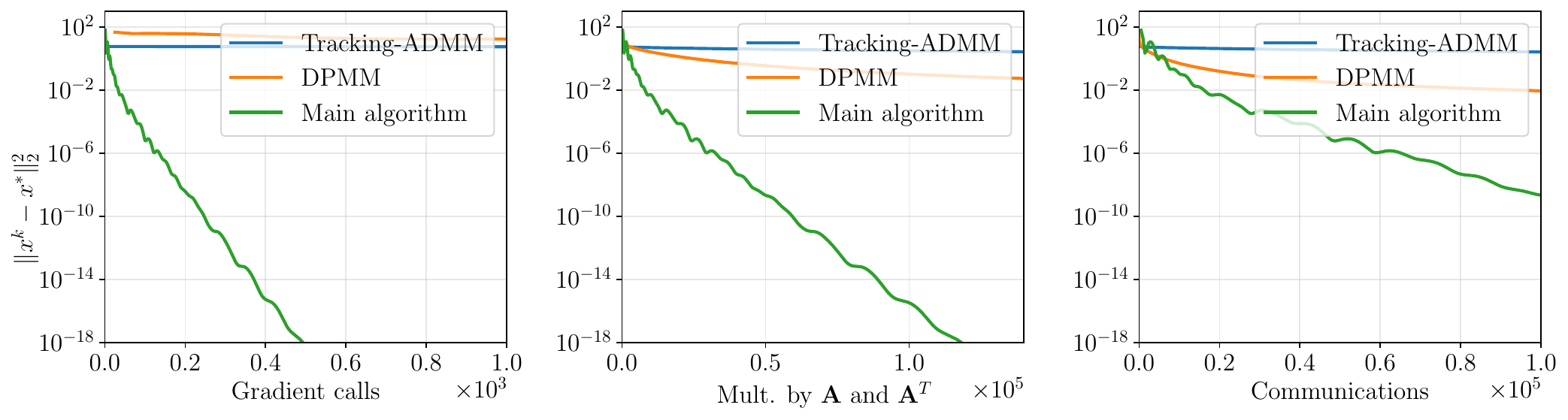}
    \vspace{-0.5cm}
    \caption{VFL, Erd\H{o}s–R\'enyi graph, $n = 7$, $m = 100$}
    \label{fig:vfl-1}
\end{figure}

For numerical simulation, we use \texttt{mushrooms} dataset from LibSVM library \cite{libsvm}. We split $m=100$ samples subset vertically between $n=7$ devices. Regularization parameter $\lambda = 10^{-2}$. The results are in the Fig.~\ref{fig:vfl-1}.

Our algorithm exhibits the best convergence rates, as evidenced by the steepest slopes.
The slopes vary for gradient calls, matrix multiplications, and communications. 
This is due to the fact that \Cref{alg:main} involves many communications per iteration, in contrast to DPMM and Tracking-ADMM, which make numerous gradient calls per iteration. 

\section{Acknowledgements}\label{sec:acknowledgements}
The work was supported by MIPT based  Center of National Technology Initiatives in the field of Artificial Intelligence for the purposes of  the "road map" of Artificial Intelligence development up to 2030 and supported by NTI Foundation (agreement No.70-2021-00207 dated 22.11.2021, identifier 000000S507521QYL0002).

\cite*{}

\bibliographystyle{unsrtnat}
\bibliography{references}

\newpage
\appendix
\section*{Appendix / supplemental material}

\section{Reduction of consensus optimization to coupled constraints}\label{sec:consensus_vs_coupled}

As mentioned in \Cref{sec:intro}, the consensus optimization problem 
\begin{align*}
  \min_{x_1,\ldots,x_n \in \sX} \sum_{i=1}^n f_i(x_i) \qst x_1 = x_2 = \ldots = x_n
\end{align*}
can be reduced to the problem with coupled constraints \eqref{eq:main}.
During the review process, we were asked whether the complexity of \Cref{alg:main} could be reduced to match that of optimal algorithms for consensus optimization \cite{scaman2017optimal}.

It turns out that, for a general-purpose first-order decentralized algorithm for problems with coupled constraints (as defined in \Cref{subsec:lower_bounds}), the communication complexity is worse by at least  a factor of $\sqrt{n}$. This conclusion is based on our complexity lower bounds in \Cref{thm:lower_bounds} and the following lower bound on $\kappa_\mA$ for any suitable matrix $\mA$.

Let $d=1$.
Consider an arbitrary horizontal-block matrix $\mA' = (\mA_1 \ldots \mA_n) \in \R^{m \times n}$ such that $\ker \mA' = \cL_1 = \Span \{\vec 1_n\}$, which corresponds to the consensus constraint.
By definition \eqref{eq:def_mu_a_L_a}, we have
\begin{equation*}
\mu_\mA = \frac1n\lambda_{\min^+}\cbr{\sum_{i=1}^n \mA_i \mA_i^\top} = \frac1n\lambda_{\min^+}\cbr{\mA' \mA'^\top} = \frac1n\sigma_{\min^+}^2(\mA'),
\end{equation*}
and
\begin{equation*}
L_\mA = \max_{i=1\ldots n}\sigma^2_{\max}(\mA_i) = \|\mA_{i^*}\|^2,
\end{equation*}
where $i^*$ denotes the index of the column achieving the maximum.
We can upper bound the minimal positive singular value using a vector $v$, whose components are all zero except for one $1$ and one $-1$, implying $\|v\|^ 2= 2$. Since the sum of its components is zero, $v$ is orthogonal to $\ker \mA'$, thus
\begin{equation*}
\sigma_{\min^+}^2(\mA') =\min_{\substack{v:~ \|v\|> 0,\\ v \bot \vec 1}} \|\mA'v\|^2 / \|v\|^2 \leq \|\mA_i -\mA_j\|^2 / 2 \leq 2 \|\mA_{i^*}\|^2.
\end{equation*}
It follows that
\begin{equation*}
\kappa_\mA = \frac{L_\mA}{\mu_\mA} \geq \frac{\|\mA_{i^*}\|^2}{\frac2n\|\mA_{i^*}\|^2} = \frac n2.
\end{equation*}
This bound is nearly tight, because we can achieve $\kappa_\mA = n-1$ by setting $\mA'$ to be the Laplacian matrix or the incidence matrix of a complete graph on $n$ nodes.

Substituting this  into the complexity lower bounds in \Cref{thm:lower_bounds}, we find that the number of communication rounds is lower bounded by 
$\Omega\cbr{\sqrt{\kappa_f}\sqrt{n}\sqrt{\kappa_\mW}\log\cbraces{\frac{1}{\eps}}}$
for any choice of $\mA$, while the optimal complexity for first order decentralized consensus optimization is  
$O\cbr{\sqrt{\kappa_f}\sqrt{\kappa_\mW}\log\cbraces{\frac{1}{\eps}}}$.

\section{Missing proofs from Section~\ref{sec:derivation_of_algorithm}}\label{app:derivation_of_algorithm}

\subsection{Proof of Lemma~\ref{lem:G}}\label{app:lem:G}
\begin{proof}
  Let $D_G(x', y'; x, y)$ denote the Bregman divergence of $G$:
	\begin{equation}\label{eq:bregman}
		D_G(x', y'; x, y) = G(x', y') - G(x, y) - \<\nabla_x G(x, y), x' -x> - \<\nabla_y G(x, y), y' - y>.
	\end{equation}

	The value of $\mu_G$ can be obtained as follows:
	\begin{align*}
		\bg_{G}(x',y';x,y)
		 & =
		\bg_F(x';x) + \frac{r}{2}\sqn{\mA(x'-x) + \gamma \pW (y'-y)}
		\\&\ageq{is due to \Cref{ass:f}}
		\frac{\mu_f}{2}\sqn{x'-x} + \frac{r}{2}\sqn{\mA(x'-x) + \gamma \pW (y'-y)}
    \\&=
    \frac{\mu_f}{2}\sqn{x'-x} + \frac{r}{2}\sqn{\mA(x'-x)} + r\<\mA(x'-x), \gamma \pW (y'-y)> 
    \\& \quad
    + \frac{r}{2}\sqn{\gamma \pW (y'-y)}
		\\&\ageq{is due to Young's inequality}
		\frac{\mu_f}{2}\sqn{x'-x} + \frac{r}{4}\sqn{\gamma \pW(y'-y)} - \frac{r}{2}\sqn{\mA(x'-x)}
    \\&\ageq{is due to \Cref{ass:A}, $y' -y \in \cL_m^\perp$,  \cref{eq:range_pW} and \cref{eq:mat_Lmu}}
		\frac{\mu_f}{2}\sqn{x'-x} + \frac{r \gamma^2 \mu_\pW}{4}\sqn{y'-y} - \frac{rL_{\mA}}{2}\sqn{x'-x}
		\\&\aeq{is due to \cref{eq:r_gamma}}
		\frac{\mu_f}{4}\sqn{x'-x} + \frac{\mu_f\gamma^2 \mu_\pW}{8 L_{\mA}}\sqn{y'-y},
    \\&\ageq{is due to \cref{eq:r_gamma}}
    \frac{\mu_f}{2} \min\braces{\frac12,  \frac{\mu_\mA + L_\mA}{4 L_\mA}}\left\|\pmat{x'-x \\ y' - y}\right\|^2,
	\end{align*}
	where \annotate.

	The value of $L_G$ can be obtained as follows:
	\begin{align*}
		\bg_{G}(x',y';x,y)
		 & =
		\bg_F(x';x) + \frac{r}{2}\sqn{\mA(x'-x) + \gamma \pW (y'-y)}
		\\&\aleq{is due to \Cref{ass:f}}
		\frac{L_f}{2}\sqn{x'-x} + \frac{r}{2}\sqn{\mA(x'-x) + \gamma \pW (y'-y)}
		\\&\aleq{is due to Young's inequality}
		\frac{L_f}{2}\sqn{x'-x} + r\sqn{\gamma \pW (y'-y)} + r\sqn{\mA(x'-x)}
    \\&\aleq{is due to \Cref{ass:A} and \cref{eq:mat_Lmu}}
		\frac{L_f}{2}\sqn{x'-x} + r \gamma^2 L_\pW \sqn{y'-y} + rL_{\mA}\sqn{x'-x}
		\\&\aeq{is due to \cref{eq:r_gamma}}
		\frac{L_f + \mu_f}{2}\sqn{x'-x} + \frac{\mu_f\gamma^2 L_\pW}{2L_{\mA}}\sqn{y'-y},
    \\&\aleq{is due to \cref{eq:r_gamma}}
    \frac{1}{2} \max\braces{L_f + \mu_f, \mu_f\frac{\mu_\mA + L_\mA}{L_\mA}\frac{L_\pW}{\mu_\pW}}\left\|\pmat{x'-x \\ y' - y}\right\|^2,
	\end{align*}
	where \annotate.
\end{proof}

\subsection{Proof of Lemma~\ref{lem:B}}\label{app:lem:B}
\begin{proof}
To obtain the formula for $L_\mB$, consider an arbitrary $z  \in (\sY)^n$:
\begin{align*}
  \sqn{\mB^\top z} &= \sqn{\mA^\top z} + \sqn{\gamma\pW z} 
  \\&\aleq{is due to \Cref{ass:A} and \cref{eq:mat_Lmu}}
  (L_\mA + \gamma^2 L_\pW) \sqn{z} 
  \\&\aeq{is due to \cref{eq:r_gamma}}
  \cbraces{L_\mA + (L_\mA + \mu_\mA)\frac{L_\pW}{\mu_\pW}}\sqn{z} ,
\end{align*}
	where \annotate.

  To derive the formula for $\mu_\mB$, first of all, note that by \cref{eq:range_pW} 
  \begin{equation}
  (\ker \mB^\top)^\bot = \range \mB = \range \mA + \range \pW = \range \mA + \cL_m^\perp.
  \end{equation}
	Let $z \in (\ker \mB^\top)^\bot = u + v$, where $u = (u_1,\ldots,u_n), v = (v_0,\ldots,v_0) \in (\sY)^n$ such that $u \in \cL_m^\perp$ and $v \in \cL_m$.

	We can show that $v_0 \in \range \mS$. In order to do that, let us show that $\<v_0, w_0> = 0$ for all $w_0 \in \ker \mS$. Let $w = (w_0,\ldots,w_0) \in \cL_m$.
	The fact that $w_0 \in \ker \mS$ and $w \in \cL_m$ implies $w \in \ker \mA \mA^\top = \ker \mA^\top$. Hence, it is easy to show that $w \in \ker \mB^\top = (\range \mB)^\perp$.
	Then, we obtain
	\begin{equation*}
		n\<v_0, w_0>
		\aeq{follows from the definition of $v$ and $w$}
		\<v,w>
		\aeq{follows from the fact that $u \in \cL_m^\perp$ and $w \in \cL_m$}
		\<u + v, w>
		=
		\<z, w>
		\aeq{follows from the fact that $z \in \range \mB$ and $w \in (\range \mB)^\perp$}
		0,
	\end{equation*}
	where \annotate. Hence, $v_0 \in \range \mS$.

	Further, we get
	\begin{align*}
		\sqn{\mB^\top z}
    & \aeq{is due to the definitions of $u$ and $v$}
		\sqn{\mA^\top (u+v)} + \sqn{\gamma\pW (u+v)}
		\\&\aeq{is due to the fact that $v \in \cL_m$}
		\sqn{\mA^\top (u+v)} + \sqn{\gamma\pW u}
    \\&\ageq{is due to \cref{eq:mat_Lmu} and \cref{eq:range_pW}}
		\sqn{\mA^\top (u+v)} + \gamma^2\mu_\pW\sqn{u}
		\\&=
		\sqn{\mA^\top u} +  \sqn{\mA^\top v} + 2\<\mA^\top u, \mA^\top v> + \gamma^2\mu_\pW\sqn{u}
		\\&\ageq{uses Young's inequality}
		-\sqn{\mA^\top u} + \frac{1}{2}\sqn{\mA^\top v} +  \gamma^2\mu_\pW\sqn{u}
    \\&\aeq{is due to the definitions of $v$ and $\mS$, and 
    $\sqn{\mA^\top v} = \left\|\pmat{\mA_1^\top v_0 \\ \vdots \\\mA_n^\top v_0}\right\|^2  =\sum_{i=1}^n \sqn{\mA_i^\top v_0} = \<v_0, \sum_{i=1}^n \mA_i \mA_i^\top v_0> =  \<v_0, n\mS v_0>$}
		-\sqn{\mA^\top u} + \frac{1}{2}\<v_0, n\mS v_0> + \gamma^2\mu_\pW\sqn{u}
		\\&\ageq{is due to \Cref{ass:A} and the definition of $v$}
		-L_{\mA}\sqn{u} + \frac{n\mu_{\mA}}{2}\sqn{v_0} + \gamma^2\mu_\pW\sqn{u}
    \\&=-L_{\mA}\sqn{u} + \frac{\mu_{\mA}}{2}\sqn{v} + \gamma^2\mu_\pW\sqn{u}
		\\&\aeq{is due to \cref{eq:r_gamma}}
		\frac{\mu_{\mA}}{2}\sqn{v} + \mu_{\mA}\sqn{u}
		\\&\ageq{is due to the definitions of $u$ and $v$}
		\frac{\mu_{\mA}}{2}\sqn{z},
	\end{align*}
	where \annotate.

\end{proof}

\subsection{Proof of Theorem~\ref{thm:main}}\label{app:thm:main}
\begin{lemma}[\cite{salim2022optimal}, Section 6.3.2]\label{lem:cheb_iter}
Let $\mM$ be a matrix with $\mu_\mM > 0$, $\mr \in \range \mM$ and $\mM v_0 = \mr$. 
Then $\cP_{\mM}(\mM^\top \mM)(v - v_0) = v - \cheb(v, \mM, \mr)$,
where $ \cheb$ is defined as \Cref{alg:cheb}.
\end{lemma}

\begin{algorithm}[H]
  \caption{$\cheb(v, \mM, \mr)$: Chebyshev iteration (\cite{gutknecht2002chebyshev}, Algorithm 4)}
	\label{alg:cheb}
	\begin{algorithmic}[1]
		\State {\bf Parameters:} $v, \mM, 
		\mr.$
    \State $n \eqdef \left\lceil\sqrt{\frac{L_\mM}{\mu_\mM}}\right\rceil$
    \State $\rho\eqdef\big(L_\mM-\mu_\mM\big)^2/16$, $\nu\eqdef(L_\mM+\mu_\mM)/2$
		\State $\delta^0 \eqdef -\nu/2$
    \State $p^0 \eqdef -{\mM}^\top ({\mM} v-\mr)/\nu$
		\State $v^1 \eqdef v + p^0$
		\For{$i=1,\ldots,n-1$}{}
		\State $\beta^{i-1}\eqdef\rho/\delta^{i-1}$
		\State $\delta^i\eqdef-(\nu+\beta^{i-1})$
		\State $p^i\eqdef\big(\mM^\top (\mM v^{i}-\mr)+\beta^{i-1}p^{i-1}\big)/\delta^i$
		\State $v^{i+1}\eqdef v^i+p^i$
		\EndFor
	\State {\bf Output:} $v^n$
	\end{algorithmic}
\end{algorithm}

\textbf{Proof of \Cref{thm:main}}
\begin{proof}
Applying \Cref{lem:poly} to $\mW$ and $\mB^\top \mB$, we derive that, due to \cref{eq:Wcheb},
it holds
\begin{equation}\label{eq:pW_Lmu}
  \lmax^2(\pW) \leq L_{\pW} = (19/15)^2,
  \quad 
  \lminp^2(\pW) \geq \mu_{\pW} = (11/15)^2,
\end{equation}
and by \cref{eq:kappa_mW} the polynomial $\cP_\mW$ has a degree of $\ceil{\sqrt\kappa_\mW}$.  
Similarly, due to \cref{eq:K}, it holds 
\begin{equation}\label{eq:K_Lmu}
  \smax^2(\mK) = \lmax(\mK^\top \mK) \leq L_{\mK} = {19/15},
  \quad
  \sminp^2(\mK) = \lminp(\mK^\top \mK) \geq \mu_\mK = {11/15}, 
\end{equation}
and since $\kappa_\mB = \frac{L_\mB}{\mu_\mB}$, 
the polynomial $\cP_{\mB}$ has a degree of $\ceil{\sqrt\kappa_\mB}$.

We implement computation of the term $\mK^\top (\mK u - \pb)$ in line~\ref{alg:apapc:line:z} of \Cref{alg:apapc} via \Cref{alg:Kcheb} by \Cref{lem:cheb_iter}:
\begin{align*}
  \mK^\top (\mK u - \pb) &= \mK^\top\mK (u - u_0) 
  = \cP_\mB(\mB^\top\mB)(u - u_0)
                       \\&= u - \cheb(u, \mB, \mb) = \Kcheb(u).
\end{align*}
Similarly, utilizing~\Cref{lem:cheb_iter}, we get 
\begin{equation}
  \pW y = \cP_{\sqrt\mW}(\mW) y = \cP_{\sqrt\mW}(\sqrt{\mW}^\top \sqrt\mW) (y - 0) = y - \cheb(y, \sqrt\mW, 0) = \mulpW(y),
\end{equation}
where $\mulpW$ is defined as \Cref{alg:mult_pW}.

Therefore, \Cref{alg:main} is equivalent to \Cref{alg:apapc}.

From \cref{eq:G_mu_L,eq:pW_Lmu}, $\frac{\mu_\mA + L_\mA}{L_\mA} \leq 2$ and $(19/11)^2 \leq 3$, we get
\begin{align}
  L_G &= \max\braces{L_f + \mu_f, \mu_f\frac{\mu_\mA + L_\mA}{L_\mA}\frac{L_\pW}{\mu_\pW}} 
  \leq \mu_f\max\braces{1 + \kappa_f, 6}, \label{eq:L_G_bound}
\\
  \mu_G &= \mu_f \min\braces{\frac12, \frac{\mu_\mA + L_\mA}{4L_\mA}} \geq \frac{\mu_f}4, \label{eq:mu_G_bound}
  \\
  \kappa_G &= \frac{L_G}{\mu_G} \leq 4\max\braces{1 + \kappa_f, 6}. \label{eq:kappa_G_bound}
\end{align}

From \cref{eq:kappa_B,eq:pW_Lmu} we get
\begin{equation}\label{eq:kappa_B_specific}
  \kappa_\mB = \frac{L_\mB}{\mu_\mB} \leq 2\cbraces{\kappa_\mA + (19/11)^2(1 + \kappa_\mA)} \leq 8\kappa_\mA + 6.
\end{equation}

From \cref{eq:K_Lmu} we obtain
\begin{equation}\label{eq:kappa_mK}
  \kappa_\mK = \frac{L_\mK}{\mu_\mK} = 19/11,
\end{equation}
and substituting \cref{eq:kappa_B_specific,eq:kappa_G_bound} to \Cref{prop:apapc}, we obtain as its direct corollary that $k = O(\sqrt{\kappa_f} \log(1/\eps))$.
Each iteration of \Cref{alg:main} require $O(1)$ computations of $\nabla F$, $O(\sqrt{\kappa_\mB}) = O(\sqrt{\kappa_\mA})$ multiplications by $\mA,\mA^\top$ and
$O(\sqrt{\kappa_\mA}\sqrt{\kappa_\mW})$ multiplications by $\mW$, which gives us the statement of \Cref{thm:main}.
The values of the parameters $\tau, \eta, \theta, \alpha$ in \Cref{thm:main} are derived from \Cref{prop:apapc} as follows.  We have $\tau = \min\braces{1, \frac12\sqrt\frac{\kappa_\mK}{\kappa_G}} = \min\braces{1, \frac12\sqrt\frac{19}{44\max\braces{1+\kappa_f,  6}}}$ due to \cref{eq:kappa_mK,eq:kappa_G_bound};
$\eta = \frac{1}{4\tau L_G} = \frac{1}{4\tau \max\braces{L_f + \mu_f, 6\mu_f}}$ due to \cref{eq:L_G_bound}; 
$\theta = \frac{15}{19\eta}$ due to \cref{eq:K_Lmu} and $\alpha = \mu_G = \frac{\mu_f}{4}$ due to \cref{eq:mu_G_bound}.
\end{proof}

\section{Proof of Theorem~\ref{thm:lower_bounds}}\label{app:lower_bounds_proof}

\subsection{Dual problem }

Let us construct the lower bound for the problem dual to the initial one. Consider primal problem with zero r.h.s. in constraints
\begin{equation}
  \begin{aligned}
   \label{prob:lower_primal}
	\min_{x_1, \ldots, x_n\in\ell_2}~ &\sum_{i=1}^n f_i(x_i) \\
	\text{s.t. } &\sum_{i=1}^n \mA_i x_i = 0.
  \end{aligned}
\end{equation}

The dual problem has the form
\begin{equation}
\begin{aligned}
  \label{prob:lower_dual}
	&\min_{x_1, \ldots, x_n\in\ell_2} \max_{z} \sbraces{\sum_{i=1}^n f_i(x_i) - \angles{z, \mA_i x_i}} 
	= \max_z \sbraces{-\max_{x_1, \ldots, x_n\in\ell_2} \sum_{i=1}^n \angles{\mA_i^\top z, x_i} - f_i(x_i)} \\
	&\quad= -\min_z \sum_{i=1}^n f_i^*(\mA_i^\top z).
\end{aligned}
\end{equation}
Introducing local copies of $w$ at each node, we get
\begin{align}\label{eq:lower_bound_dual_problem}
\begin{aligned}
	\min_{z_1, \ldots, z_n}~ &\sum_{i=1}^n g_i(z_i) := \sum_{i=1}^n f_i^* (\mA_i^\top z_i) \\
	&\text{s.t. } \mW z = 0.
\end{aligned}
\end{align}

\subsection{Example  graph}

We follow the principle of lower bounds construction introduced in \cite{kovalev2021lower} and take the example graph from \cite{scaman2017optimal}. Let the functions held by the nodes be organized into a path graph with $n$ vertices, where $n$ is divisible by $3$. The nodes of graph $\cG = (\cV, \cE)$ are divided into three groups $\cV_1 = \braces{1, \ldots, n/3}, \cV_2 = \braces{n/3 + 1, \ldots, 2n/3}, \cV_3 = \braces{2n/3 + 1, \ldots, n}$ of $n/3$ vertices each.

Now we recall the construction from \cite{scaman2017optimal}. Maximum and minimum eigenvalues of a path graph have form $\lambda_{\max}(W) = 2\cbraces{1 + \cos\frac{\pi}{n}},~ \lambda_{\min^+}(W) = 2\cbraces{1 - \cos\frac{\pi}{n}}$. Let $\beta_n = \frac{1 + \cos\cbraces{\frac{\pi}{n}}}{1 - \cos\cbraces{\frac{\pi}{n}}}$. Since $\beta_n\overset{n\to\infty}{\rightarrow} +\infty$, there exists $n = 3m\geq 3$ such that $\beta_n\leq \kappa_\mW < \beta_{n+3}$. For this $n$, introduce edge weights $w_{i, i+1} = 1 - a\mathbb{I}\braces{i=1}$, take the corresponding weighed Laplacian $W_a$ and denote its condition number $\kappa(W_a)$. If $a = 1$, the network is disconnected and therefore $\kappa(W_a) = \infty$. If $a = 0$, we have $\kappa(W_a) = \beta_n$.  By continuity of Laplacian spectra we obtain that for some $a \in [0, 1)$ it holds $\kappa(W_a) = \kappa_\mW$. Note that $\pi/(n+3)\in[0, \pi/3]$ for $x\in[0, \pi/3]$ it holds $1 - \cos x\geq x^2/4$. We have
\begin{align}\label{eq:n_chi_relation}
	\kappa_\mW\leq \beta_{n+3} = \frac{1 + \cos\frac{\pi}{n+3}}{1 - \cos\frac{\pi}{n+3}}
	\leq \frac{72(n+3)^2}{\pi^2}\leq \frac{288 n^2}{\pi^2}\leq 32n^2~~~\Rightarrow~~~ \sqrt{\kappa_\mW}\leq 4\sqrt 2n = O(n).
\end{align}

\subsection{Example functions}
We let $e_1 = (1~ 0~ \ldots~ 0)^\top$ denote the first coordinate vector and define functions
\begin{align}\label{eq:lower_bound_functions_f}
	f_i(p, t) = \frac{\mu_f}{2} \left\|p + \frac{\sqrt{\hat L_\mA}}{2\mu_f}e_1\right\|^2 + \frac{L_f}{2} \sqn{t}.
\end{align}
We immediately note that each $f_i$ is $L_f$-smooth and $\mu_f$-strongly convex. The first term has the form $\frac {\mu_f}2 \sqn{p + c e_1}$, and its convex conjugate is
$\cbr{\frac{\mu_f}2 \sqn{p + c e_1}}^* = \max_p \braces{\<u, p> - \frac{\mu_f}{2}\sqn{p + c e_1}}$. 
The gradient by $p$ is $u - \mu_f(p + c e_1) = 0$, thus 
$p = \frac{1}{\mu_f} u - c e_1$ and 
$\max_p \braces{\<u, p> - \frac{\mu_f}{2}\sqn{p + c e_1}} = \frac{\sqn u}{2 \mu_f} - cu_1$.
Correspondingly,
\begin{align*}
	f_i^*(u, v) = \frac{1}{2\mu_f} \sqn{u} + \frac{1}{2L_f} \sqn{v} - \frac{\sqrt{\hat L_\mA}}{2\mu_f}u_1.
\end{align*}
To define matrices $\mA_i$, we first introduce
\begin{align*}
	\mE_1 = \begin{pmatrix}
		1 & 0 & 0 & 0 & 0 & \ldots & \\
		0 & 1 & -1 & 0 & 0 & \ldots & \\
		0 & 0 & 0 & 0 & 0 & \ldots & \\
		0 & 0 & 0 & 1 & -1 & \ldots & \\
		\vdots & \vdots & \vdots & \vdots & \vdots & \ddots \\		
	\end{pmatrix},~~~~
	\mE_2 = \begin{pmatrix}
		1 & -1 & 0 & 0 & 0 & \ldots & \\
		0 & 0 & 0 & 0 & 0 & \ldots & \\
		0 & 0 & 1 & -1 & 0 & \ldots & \\
		0 & 0 & 0 & 0 & 0 & \ldots & \\
		\vdots & \vdots & \vdots & \vdots & \vdots & \ddots & \\		
	\end{pmatrix}.
\end{align*}
Let $\hat L_\mA = \frac{1}{2} L_\mA - \frac{3}{4} \mu_\mA,~ \hat\mu_\mA = \frac{3}{2}\mu_\mA$ and introduce
\begin{align*}
	\mA_i = \begin{cases}
		[\sqrt{\hat L_\mA} \mE_1^\top ~~ \sqrt{\hat\mu_\mA} \mI], &i\in\cV_1 \\
		[~~~\mathbf{0}~~~~~~~~~~~~~~~~~  \mathbf{0}~~~],  &i\in\cV_2 \\
		[\sqrt{\hat L_\mA} \mE_2^\top ~~ \sqrt{\hat\mu_\mA} \mI], &i\in\cV_3
	\end{cases}
\end{align*}
Let us make sure that the choice of $\mA_i$ guarantees constants $L_\mA, \mu_\mA$ from \eqref{eq:def_mu_a_L_a}.
\begin{align*}
	\max_i \lambda_{\max} (\mA_i \mA_i^\top ) &= \lambda_{\max} \cbraces{\hat L_\mA \mE_1^\top \mE_1 + \hat \mu_\mA \mI} = 2\hat L_\mA + \hat\mu_\mA = L_\mA, \\
	\lambda_{\min}^+ \cbraces{\frac{1}{n} \sum_{i=1}^n \mA_i \mA_i^\top} 
	&= \lambda_{\min}^+ \cbraces{\frac{1}{3}(\hat L_\mA \mE_1^\top \mE_1 + \hat \mu_\mA\mI) + \frac{1}{3}(\hat L_\mA \mE_2^\top \mE_2 + \hat \mu_\mA\mI)} 
	\\
	&= \frac{2}{3} \hat\mu_\mA = \mu_\mA.
\end{align*}
Let $\widetilde\mM = \begin{pmatrix} ~~~1 & -1 \\ -1 & ~~~1 \end{pmatrix}$ and
\begin{align*}
	\mM_1 &= \mE_1^\top \mE_1 = 
	\begin{pmatrix}
		1 & 0 & 0 & \ldots \\
		0 & \widetilde\mM & 0 & \ldots \\
		0 & 0 & \widetilde\mM & \ldots \\
		\vdots & \vdots & \vdots & \ddots
	\end{pmatrix},~~~~
	\mM_2 = \mE_2^\top \mE_2 = 
	\begin{pmatrix}
		\widetilde\mM & 0 & 0 & \ldots \\
		0 & \widetilde\mM & 0 & \ldots \\
		0 & 0 & \widetilde\mM & \ldots \\
		\vdots & \vdots& \vdots& \ddots
	\end{pmatrix}
\end{align*}
The dual functions take the form
\begin{align}
	g_i(z) = f_i^*(\mA_i^\top z) &= 
	\begin{cases}
		\frac{1}{2\mu_f}\sqn{\sqrt{\hat L_\mA} \mE_1 z} + \frac{1}{2L_f}\sqn{\sqrt{\hat\mu_\mA} z} - \frac{\hat L_\mA}{2\mu_f} z_1, &i\in\cV_1 \\
		0, &i\in\cV_2 \\
		\frac{1}{2\mu_f}\sqn{\sqrt{\hat L_\mA} \mE_2 z} + \frac{1}{2L_f}\sqn{\sqrt{\hat\mu_\mA} z} - \frac{\hat L_\mA}{2\mu_f} z_1, &i\in\cV_3
	\end{cases}  \nonumber \\
	&=
	\begin{cases}
		\frac{1}{2}z^\top \cbraces{\frac{\hat L_\mA}{\mu_f} \mM_1 + \frac{\hat\mu_\mA}{L_f} \mI} z - \frac{\hat L_\mA}{2\mu_f} z_1, &i\in\cV_1 \\
		0, &i\in\cV_2 \\
		\frac{1}{2}z^\top \cbraces{\frac{\hat L_\mA}{\mu_f} \mM_2 + \frac{\hat\mu_\mA}{L_f} \mI} z - \frac{\hat L_\mA}{2\mu_f} z_1, &i\in\cV_3
	\end{cases} \label{eq:g_functions_example}
\end{align}
Therefore, we have
\begin{align}\label{eq:lower_bound_sum_g_i}
\begin{aligned}
	\sum_{i=1}^n g_i(z) &= \frac{n}{3} \sbraces{\frac{\hat L_\mA}{2\mu_f} z^\top (\mM_1 + \mM_2) z + \frac{\hat\mu_\mA}{L_f} z^\top z - \frac{\hat L_\mA}{\mu_f} z_1} \\
	&= \frac{n}{3} \frac{\hat L_\mA}{\mu_f} \sbraces{\frac{1}{2} z^\top \mM z - z_1 + \frac{\hat \mu_\mA \mu_f}{\hat L_\mA L_f} z^\top z},
\end{aligned}
\end{align}
where
\begin{align*}
	\mM = \mM_1+ \mM_2 = \begin{pmatrix}
		~~~2 & -1 & ~~~0 & ~~~0 & 0 & \ldots \\
		-1 & ~~~2 & -1 & ~~~0 & 0 & \ldots \\
		~~~0 & -1 & ~~~2 & -1 & 0 & \ldots \\
		\vdots & \vdots & \vdots & \vdots & \vdots & \ddots
	\end{pmatrix}.
\end{align*}
Now we formulate the lower complexity bounds for $\sum_{i=1}^n g_i(z)$, where $g_i(z)$ are defined in \eqref{eq:g_functions_example}.

\subsection{Deriving the lower bound}
\begin{lemma}
	Function $\sum_{i=1}^n g_i(z)$ attains its minimum at $z^* = \braces{\rho^k}_{k=1}^\infty$, where
	\begin{align*}
		\rho = \frac{\sqrt{\frac{2L_\mA L_f}{\mu_\mA \mu_f} + 1} - 1}{\sqrt{\frac{2L_\mA L_f}{\mu_\mA \mu_f} + 1} + 1}.
	\end{align*}
\end{lemma}
\begin{proof}
In the lower bound example in (Lemma~1 in Appendix~C) \cite{kovalev2021lower} it was shown that function
\begin{align*}
	h(z) = \frac{1}{2} z^\top \mM z + \frac{3\mu}{L - \mu} z^\top z - z_1
\end{align*}
attains its minimum at $z_k^* = \rho^k$, where
\begin{align*}
	\rho = \frac{\sqrt{\frac{2L}{3\mu} + \frac{1}{3}} - 1}{\sqrt{\frac{2L}{3\mu} + \frac{1}{3}} + 1}.
\end{align*}
Let us deduce the expression for $L / \mu$ in terms of $L_\mA, L_f, \mu_\mA, \mu_f$. We enforce $h(z) = \sum_{i=1}^n g_i(z)$ and set
\begin{align*}
	\frac{3\mu}{L - \mu} = \frac{\mu_\mA \mu_f}{L_\mA L_f}\Rightarrow \frac{L}{\mu} = 1 + 3\frac{L_\mA L_f}{\mu_\mA \mu_f}.
\end{align*}
Therefore, $h(z)$ attains its minimum at $z^k = \rho^k$, where
\begin{align*}
	\rho = \frac{\sqrt{\frac{2L_\mA L_f}{\mu_\mA \mu_f} + 1} - 1}{\sqrt{\frac{2L_\mA L_f}{\mu_\mA \mu_f} + 1} + 1}.
\end{align*}
\end{proof}
Let us first show the lower bound on the number of communications.
Without loss of generality we can assume that the initial point chosen by a first-order algorithm is $x^0 = 0$ (otherwise we can  shift the variables accordingly), thus $\cM_i(0) = \{0\},~ \cH_i(0) = \braces{0}$.
\begin{lemma}
	Let $s_i(k)$ denote the maximum index of a nonzero component of vector blocks $p, t$ held by $i$-th node at step $k$ in its local memory $\cM_i(k)$ and vector $z$ held by $i$-th node in its local memory $\cH_i(k)$, i.e.
	\begin{align*}
		s_i(k) = 
		\begin{cases}
			0, &\hspace{-1.0cm}\text{if } \cM_i(k)\subseteq\braces{0} \text{ and } \cH_i(k)\subseteq\braces{0} \\
			\min\Big\{s\in \braces{1, 2, \ldots}: \cH_i(k)\subseteq\Span\braces{e_1, \ldots, e_s} \\ \qquad~ \text{ and } \cM_i(k)\subseteq \Span\braces{e_1, \ldots, e_s} \times\braces{e_1, \ldots, e_s}\Big\}, &\text{else}.
		\end{cases}
	\end{align*}
	Let $k_q$ denote the number of algorithm step by which exactly $q$ communication steps have been performed, where $q\geq 0$. For any $k\in\braces{1, \ldots, k_q}$ we have
	\begin{align}\label{eq:nonzero_components_bound}
		\max_i s_i(k)\leq 2 + \left\lfloor \frac{q}{\frac{n}{3}+ 1} \right\rfloor
	\end{align}
\end{lemma}
\begin{proof}
	If the method performs a multiplication by $\mA_i$, it transfers the information from $\cM_i(k)$ to $\cH_i(k)$. If multiplication by $\mA_i^\top$ is performed, the information is transferred in the opposite direction, i.e. from $\cH_i(k)$ to $\cM_i(k)$. If the method performs a matrix multiplication step (either by $\mA$ or $\mA^\top$), then from the structure of $\mA_i$ we obtain
	\begin{align*}
		s_i(k+1)\leq s_i(k) + \begin{cases} 1 - (s_i(k) \mod 2), &i\in\cV_1 \\ 0, &i\in\cV_2 \\ (s_i(k)\mod 2), &i\in\cV_3 \end{cases}
	\end{align*}
	We will prove that\\
	1. For $q = 2\ell(n/3 + 1),~ \ell\in\braces{0, 1, \ldots}$ we have
	\begin{align}\label{eq:info_flow_even}
		s_i(k_q)\leq
		\begin{cases}
			1 + 2\ell,~~ i\in\cV_1 \\
			1 + 2\ell,~~ i\in\cV_2 \\
			2 + 2\ell,~~ i\in\cV_3
		\end{cases}
	\end{align}
	2. For $q = (2\ell + 1)(n/3 + 1),~ \ell\in\braces{0, 1, \ldots}$ we have
	\begin{align}\label{eq:info_flow_odd}
		s_i(k_q)\leq
		\begin{cases}
			2 + (2\ell + 1),~~ i\in\cV_1 \\
			1 + (2\ell + 1),~~ i\in\cV_2 \\
			1 + (2\ell + 1),~~ i\in\cV_3
		\end{cases}
	\end{align}
	The proof follows by induction.
	
	\noindent\textbf{Induction basis}. Let $q = 0$. From definitions of $f_i$ and $\mA_i$ it follows that
	\begin{align*}
			s_i(k_0) \leq 
			\begin{cases}
			1,~~ i\in\cV_1 \\
			0,~~ i\in\cV_2 \\
			2,~~ i\in\cV_3
		\end{cases}
	\end{align*}
	Therefore, for $q = 0$ our statement holds.
	
	\noindent\textbf{Induction step for $q = (2\ell + 1)(n/3 + 1)$}. Consider $q_- = q - n/3 = 2\ell(n/3 + 1)$. From \eqref{eq:info_flow_even} we have that for the spread of nonzero components from $\cV_3$ to $\cV_1$ it requires $n/3$ communication rounds to reach node $n/3 + 1$. After one more communication round, the information reaches node $n/3$.
	
	\noindent\textbf{Induction step for $q = 2\ell (n/3 + 1)$}. The proof follows by the same argument as for $q = (2\ell + 1)(n/3 + 1)$.

	\vspace{0.5cm} We just proved the statement of lemma, i.e. relation \eqref{eq:nonzero_components_bound}, for $q$ divisible by $(n/3 + 1)$. Between such checkpoints, the information (i.e. the number of nonzero components) traverses nodes of $\cV_2$ and therefore $\max_i s_i(k)$ stays unchanged. Thus the statement of the lemma is proven. 
\end{proof}

Now we estimate the distance to optimum.
Due to the strong duality, the solution of problem \eqref{prob:lower_primal} can be obtained from the solution of its dual $\eqref{prob:lower_dual}$ as
\begin{equation*}
x_i^*(z^*) = \pmat{p_i \\ t_i}(\mA_i^\top z^*) = 
\pmat{\frac{\sqrt{\hat L_\mA}}{\mu_f} (E_{(i)} z^* - e_1) 
\\
\frac{\sqrt{\hat\mu_\mA}}{L_f}z^*
}. 
\end{equation*}
Therefore, denoting $\alpha = \frac{\hat\mu_\mA}{L_f^2}$,
\begin{align*}
  \sqn{x_i^k - x_i^*} 
  &\geq
  \sqn{t_i^k - t_i^*} 
  \geq
  \sum_{\ell = s_i(k) + 1}^\infty (t_{i,\ell}^*)^2
  =
  \alpha \sum_{\ell = s_i(k) + 1}^\infty (z_\ell^*)^2 
  = 
  \alpha\sum_{\ell = s_i(k) + 1}^\infty \rho^{2\ell} 
  \\ &= 
  \alpha\frac{\rho^{2s_i(k) + 2}}{1 - \rho^2} 
  = 
  \alpha\frac{\rho^{6 + 2\left\lfloor\frac{q}{n/3 + 1}\right\rfloor}}{1 - \rho^2} 
  \ageq{holds since $n/3\geq 1$} 
  \alpha\frac{\rho^{6 + \frac{2q}{2n/3}}}{1 - \rho^2} 
  =
  \alpha\frac{\rho^6}{1 - \rho^2}\cdot \rho^{\frac{3q}{n}} 
  \ageq{holds due to \eqref{eq:n_chi_relation}} 
  \alpha\frac{\rho^6}{1 - \rho^2}\cdot \rho^{\frac{12q\sqrt 2}{\sqrt{\kappa_\mW}}},
\end{align*}
where \annotate.

Following \cite{kovalev2021lower}, we obtain that 
\begin{align*}
	\rho\geq \max\cbraces{0, 1 - 3\sqrt{\frac{2\mu_\mA \mu_f}{L_\mA L_f}}}.
\end{align*}
Therefore,
\begin{align*}
	\norm{x_i^k - x_i^*}_2^2\geq \alpha\frac{\rho^6}{1 - \rho^2} \cbraces{\max\cbraces{0, 1 - 3\sqrt{\frac{2\mu_\mA \mu_f}{L_\mA L_f}}}}^{\frac{3q}{\sqrt{\kappa_\mW}}}.
\end{align*}
It follows that the number of communications is lower bounded as
\begin{align}\label{eq:lower_bound_communications}
	N_{\mW}\geq \Omega\cbraces{\sqrt{\kappa_\mW}\sqrt{\frac{L_\mA L_f}{\mu_\mA \mu_f}} \log\cbraces{\frac{1}{\eps}}},
\end{align}
and the number of matrix multiplications at each node  is lower bounded as
\begin{align}\label{eq:lower_bound_matrix_communications}
	N_{\mA}\geq \Omega\cbraces{\sqrt{\frac{L_\mA L_f}{\mu_\mA \mu_f}} \log\cbraces{\frac{1}{\eps}}}.
\end{align}

\subsection{Lower bound on the number of gradient computations}

To get the lower bound on local gradient calls, let us consider a problem
\begin{align}\label{prob:lower_united_function}
\begin{aligned}
	\min_{\substack{x_1, \ldots x_n\in \R^d \\ u_1, \ldots, u_n \in \R^d}}~ &\sum_{i=1}^n f_i(x_i) + \sum_{i=1}^n v_i(u_i) \\
	\text{s.t. }~ &\sum_{i=1}^n \mA_i x_i = 0
\end{aligned}
\end{align}
where all $f_i(x)$ are defined in \eqref{eq:lower_bound_functions_f} and all $v_i(u_i)$ are similar and defined as
\begin{align*}
	v_i(u_i) = \frac{L_f}{8} u_i^\top \mM u_i + \frac{\mu_f}{2} u^\top u - \frac{L_f}{4} u_1.
\end{align*}
First, note that each $v_i(u_i)$ is $L_f$-smooth and $\mu_f$-strongly convex, i.e. problem \eqref{prob:lower_united_function} satisfies the assumptions of Theorem~\ref{thm:lower_bounds}.

Problem \eqref{prob:lower_united_function} falls into two independent parts. The first part is minimization of $\sum_{i=1}^n f_i(x_i)$ subject to constraints and is identical to \eqref{prob:lower_primal}, while the second part is minimization of $\sum_{i=1}^n v_i(u_i)$ without constraints. Therefore, the lower bounds on the number of communications and number of matrix multiplications for problem \eqref{prob:lower_united_function} are inherited from lower bounds for problem \eqref{prob:lower_primal} and are the given by $N_\mW$ in \eqref{eq:lower_bound_communications} and $N_\mA$ in \eqref{eq:lower_bound_matrix_communications}. It remains to lower bound the number of oracle calls for minimization of $\sum_{i=1}^n v_i(u_i)$.

The structure of $v_i(u_i)$ is the same as the structure of $\sum_{j=1}^n g_j(z)$ defined in \eqref{eq:lower_bound_sum_g_i}. Therefore, the minimum of each $v_i(u_i)$ is attained at the same point $u_i^* = u^*$, and the $k$-th component of $u^*$ is given by $(u^*)_k = \nu^k$, where 
\begin{align*}
	\nu = \frac{\sqrt{\frac{L_f}{\mu_f} + \frac23} - 1}{\sqrt{\frac{L_f}{\mu_f} + \frac23} + 1}.
\end{align*}

Since there is no communication constraint on $u_i$, each node runs optimization process individually. Following the same arguments as for function $\sum_{j=1}^n g_j(z)$, we get the lower bound on the number of oracle calls
\begin{align*}
	N_f\geq \Omega\cbraces{\sqrt{\frac{L_f}{\mu_f}} \log\cbraces{\frac{1}{\eps}}}.
\end{align*}



\end{document}

%% file: cmd.tex
\DeclarePairedDelimiter\ceil{\lceil}{\rceil}

\newtheorem{theorem}{Theorem}
\newtheorem{lemma}{Lemma}

\newtheorem{assumption}{Assumption}
\newtheorem{proposition}{Proposition}

\crefname{assumption}{assumption}{assumptions}

\DeclareMathOperator{\range}{\mathrm{range}}
\DeclareMathOperator{\col}{\mathrm{col}}

\newcommand{\qst}{\quad \text{s.t.} \quad}

\newcommand{\Span}{\mathrm{Span}}

\newcommand{\lmax}{\lambda_{\max}}

\newcommand{\smax}{\sigma_{\max}}
\newcommand{\sminp}{\sigma_{\min^+}}
\newcommand{\lminp}{\lambda_{\min^+}}

\newcommand{\bg}{\mathrm{D}}

\newcommand{\R}{\mathbb{R}}

\def\<#1,#2>{\langle #1,#2\rangle}

\newcommand{\norm}[1]{\|#1\|}
\newcommand{\sqn}[1]{\norm{#1}^2}
\newcommand{\Norm}[1]{\left\|#1\right\|}
\newcommand{\sqN}[1]{\Norm{#1}^2}

\newcommand{\diag}[1]{\mathrm{diag}\left(#1\right)}

\newcommand{\sX}{\R^d}
\newcommand{\sY}{\R^m}
\newcommand{\sU}{{\mathcal U}}

\newcommand{\rng}[2]{\{#1,\ldots,#2\}}

\newcommand{\cM}{\mathcal{M}}
\newcommand{\cH}{\mathcal{H}}

\newcommand{\cG}{\mathcal{G}}

\newcommand{\cV}{\mathcal{V}}

\newcommand{\cU}{\mathcal{U}}
\newcommand{\cE}{\mathcal{E}}
\newcommand{\cL}{\mathcal{L}}

\newcommand{\cP}{\mathcal{P}}
\newcommand{\mA}{\mathbf{A}}
\newcommand{\mB}{\mathbf{B}}
\newcommand{\mK}{\mathbf{K}}
\newcommand{\mC}{\mathbf{C}}
\newcommand{\mF}{\mathbf{F}}
\newcommand{\mW}{\mathbf{W}}
\newcommand{\pW}{\mathbf{W'}}
\newcommand{\mS}{\mathbf{S}}

\newcommand{\mI}{\mathbf{I}}
\newcommand{\mE}{\mathbf{E}}

\newcommand{\mb}{\mathbf{b}}
\newcommand{\mr}{\mathbf{r}}
\newcommand{\pb}{\mathbf{b'}}

\newcommand{\mM}{\mathbf{M}}

\newcommand{\eqdef}{\coloneqq}

\newcommand{\mulpW}{\mathrm{\mathbf{mul}\pW}}
\newcommand{\cheb}{\mathrm{\mathbf{Chebyshev}}}
\newcommand{\Kcheb}{\mathrm{\mathbf{K\_Chebyshev}}}
\newcommand{\gradG}{\mathrm{\mathbf{grad\_G}}}

\newcommand{\pmat}[1]{\begin{pmatrix}#1\end{pmatrix}}

\newcommand{\angles}[1]{\left\langle #1 \right\rangle}
\newcommand{\cbraces}[1]{\left( #1 \right)}
\newcommand{\cbr}[1]{\left( #1 \right)}
\newcommand{\sbraces}[1]{\left[ #1 \right]}
\newcommand{\braces}[1]{\left\{ #1 \right\}}
\newcommand{\eps}{\varepsilon}

\newcommand{\bd}{\mathbf{d}}

\newcounter{annotatecount}
\newcounter{annotateidx}
\newcounter{annotatejdx}
\newcounter{annotatelabelcount}

\newcommand{\atran}[2]{\stepcounter{annotatecount}\overset{(\alph{annotatecount})}{#1}\csgdef{annotatedescription\theannotatecount}{#2}}
\newcommand{\aeq}[1]{\atran{=}{#1}}
\newcommand{\aleq}[1]{\atran{\leq}{#1}}
\newcommand{\ageq}[1]{\atran{\geq}{#1}}
\newcommand{\alra}[1]{\atran{\Leftrightarrow}{#1}}

\newcommand{\annotateinitused}{\setcounter{annotateidx}{0}\whileboolexpr{test{\ifnumless{\theannotateidx}{\theannotatecount}}}{\stepcounter{annotateidx}\csgdef{aused\theannotateidx}{0}}}
\newcommand{\annotategetlabels}{\setcounter{annotatejdx}{0}\setcounter{annotatelabelcount}{0}\whileboolexpr{test{\ifnumless{\theannotatejdx}{\theannotatecount}}}{\stepcounter{annotatejdx}\ifcsequal{annotatedescription\theannotateidx}{annotatedescription\theannotatejdx}{\csgdef{aused\theannotatejdx}{1}\stepcounter{annotatelabelcount}\csedef{annotatelabel\theannotatelabelcount}{(\alph{annotatejdx})}}{}}}
\newcommand{\annotateprintlabels}{\setcounter{annotatejdx}{0}\whileboolexpr{test{\ifnumless{\theannotatejdx}{\theannotatelabelcount}}}{\stepcounter{annotatejdx}\ifnumequal{\theannotatejdx}{\theannotatelabelcount}{\ifnumequal{\theannotatejdx}{1}{}{~and~}}{}\csuse{annotatelabel\theannotatejdx}\ifnumless{\theannotatejdx}{\theannotatelabelcount}{\ifnumless{\theannotatejdx+1}{\theannotatelabelcount}{,~}{}}{}}}
\newcommand{\annotate}{\annotateinitused\setcounter{annotateidx}{0}\whileboolexpr{test{\ifnumless{\theannotateidx}{\theannotatecount}}}{\stepcounter{annotateidx}\ifcsstring{aused\theannotateidx}{0}{\ifnumequal{\theannotateidx}{1}{}{;~}\annotategetlabels\annotateprintlabels~\csuse{annotatedescription\theannotateidx}}{}}\setcounter{annotatecount}{0}}

\usepackage[textsize=tiny]{todonotes}

%% file: references.bib
@article{dvinskikh2019decentralized,
	title={Decentralized and parallel primal and dual accelerated methods for stochastic convex programming problems},
	author={Dvinskikh, Darina and Gasnikov, Alexander},
	journal={Journal of Inverse and Ill-posed Problems},
	volume={29},
	number={3},
	pages={385--405},
	year={2021}
}

@inproceedings{lian2017can,
	title={Can decentralized algorithms outperform centralized algorithms? a case study for decentralized parallel stochastic gradient descent},
	author={Lian, Xiangru and Zhang, Ce and Zhang, Huan and Hsieh, Cho-Jui and Zhang, Wei and Liu, Ji},
	booktitle={Advances in Neural Information Processing Systems},
	pages={5330--5340},
	year={2017}
}

@inproceedings{ram2009distributed,
	title={Distributed non-autonomous power control through distributed convex optimization},
	author={Ram, Sundhar Srinivasan and Veeravalli, Venugopal V and Nedic, Angelia},
	booktitle={IEEE INFOCOM 2009},
	pages={3001--3005},
	year={2009},
	organization={IEEE}
}

@article{boyd2011distributed,
	author = {Boyd, Stephen and Parikh, Neal and Chu, Eric and Peleato, Borja and Eckstein, Jonathan},
	title = {Distributed Optimization and Statistical Learning via the Alternating Direction Method of Multipliers},
	journal = {Found. Trends Mach. Learn.},
	issue_date = {January 2011},
	volume = {3},
	number = {1},
	month = jan,
	year = {2011},
	issn = {1935-8237},
	pages = {1--122},
	numpages = {122},
	url = {http://dx.doi.org/10.1561/2200000016},
	doi = {10.1561/2200000016},
	acmid = {2185816},
	publisher = {Now Publishers Inc.},
	address = {Hanover, MA, USA},
}

@book{nesterov2004introduction,
	author = {Nesterov, Yurii},
	title = {Introductory Lectures on Convex Optimization: a basic course},
	publisher = {Kluwer Academic Publishers, Massachusetts},
	year = {2004}
}

@inproceedings{scaman2017optimal,
	title={Optimal algorithms for smooth and strongly convex distributed optimization in networks},
	author={Scaman, Kevin and Bach, Francis and Bubeck, S{\'e}bastien and Lee, Yin Tat and Massouli{\'e}, Laurent},
	booktitle={Proceedings of the 34th International Conference on Machine Learning-Volume 70},
	pages={3027--3036},
	year={2017},
	organization={JMLR. org}
}

@book{bertsekas1989parallel,
	title={Parallel and distributed computation: numerical methods},
	author={Bertsekas, Dimitri P and Tsitsiklis, John N},
	volume={23},
	year={1989},
	publisher={Prentice hall Englewood Cliffs, NJ}
}

@article{kovalev2020optimal,
	title={Optimal and practical algorithms for smooth and strongly convex decentralized optimization},
	author={Kovalev, Dmitry and Salim, Adil and Richt{\'a}rik, Peter},
	journal={Advances in Neural Information Processing Systems},
	volume={33},
	year={2020}
}

@article{nedic2020distributed,
	title={Distributed Gradient Methods for Convex Machine Learning Problems in Networks: Distributed Optimization},
	author={Nedi{\'c}, Angelia},
	journal={IEEE Signal Processing Magazine},
	volume={37},
	number={3},
	pages={92--101},
	year={2020},
	publisher={IEEE}
}

@article{nedic2017achieving,
	title={Achieving geometric convergence for distributed optimization over time-varying graphs},
	author={Nedic, Angelia and Olshevsky, Alex and Shi, Wei},
	journal={SIAM Journal on Optimization},
	volume={27},
	number={4},
	pages={2597--2633},
	year={2017},
	publisher={SIAM}
}

@article{nedic2009distributed,
	title={Distributed subgradient methods for multi-agent optimization},
	author={Nedi{\'c}, Angelia and Ozdaglar, Asuman},
	journal={IEEE Transactions on Automatic Control},
	volume={54},
	number={1},
	pages={48--61},
	year={2009},
	publisher={IEEE}
}

@article{olshevsky2009convergence,
	title={Convergence speed in distributed consensus and averaging},
	author={Olshevsky, Alex and Tsitsiklis, John N},
	journal={SIAM Journal on Control and Optimization},
	volume={48},
	number={1},
	pages={33--55},
	year={2009},
	publisher={SIAM}
}

@article{boyd2006randomized,
	title={Randomized gossip algorithms},
	author={Boyd, Stephen and Ghosh, Arpita and Prabhakar, Balaji and Shah, Devavrat},
	journal={IEEE transactions on information theory},
	volume={52},
	number={6},
	pages={2508--2530},
	year={2006},
	publisher={IEEE}
}

@techreport{tsitsiklis1984problems,
	title={Problems in decentralized decision making and computation.},
	author={Tsitsiklis, John Nikolas},
	year={1984},
	institution={Massachusetts Inst of Tech Cambridge Lab for Information and Decision Systems}
}

@article{li2021accelerated,
	title={Accelerated gradient tracking over time-varying graphs for decentralized optimization},
	author={Li, Huan and Lin, Zhouchen},
	journal={arXiv preprint arXiv:2104.02596},
	year={2021}
}

@article{kovalev2021lower,
	title={Lower bounds and optimal algorithms for smooth and strongly convex decentralized optimization over time-varying networks},
	author={Kovalev, Dmitry and Gasanov, Elnur and Gasnikov, Alexander and Richtarik, Peter},
	journal={Advances in Neural Information Processing Systems},
	volume={34},
	year={2021}
}

@article{rogozin2021decentralized,
	title={Decentralized distributed optimization for saddle point problems},
	author={Rogozin, Alexander and Beznosikov, Alexander and Dvinskikh, Darina and Kovalev, Dmitry and Dvurechensky, Pavel and Gasnikov, Alexander},
	journal={arXiv preprint arXiv:2102.07758},
	year={2021}
}

@article{necoara2011parallel,
	title={Parallel and distributed optimization methods for estimation and control in networks},
	author={Necoara, Ion and Nedelcu, Valentin and Dumitrache, Ioan},
	journal={Journal of Process Control},
	volume={21},
	number={5},
	pages={756--766},
	year={2011},
	publisher={Elsevier}
}

@article{necoara2014distributed,
	title={Distributed dual gradient methods and error bound conditions},
	author={Necoara, Ion and Nedelcu, Valentin},
	journal={arXiv preprint arXiv:1401.4398},
	year={2014}
}

@article{necoara2015linear,
	title={On linear convergence of a distributed dual gradient algorithm for linearly constrained separable convex problems},
	author={Necoara, Ion and Nedelcu, Valentin},
	journal={Automatica},
	volume={55},
	pages={209--216},
	year={2015},
	publisher={Elsevier}
}

@incollection{gorbunov2022recent,
	title={Recent theoretical advances in decentralized distributed convex optimization},
	author={Gorbunov, Eduard and Rogozin, Alexander and Beznosikov, Aleksandr and Dvinskikh, Darina and Gasnikov, Alexander},
	booktitle={High-Dimensional Optimization and Probability},
	pages={253--325},
	year={2022},
	publisher={Springer}
}

@article{arrow1954existence,
  title={Existence of an equilibrium for a competitive economy},
  author={Arrow, Kenneth J and Debreu, Gerard},
  journal={Econometrica: Journal of the Econometric Society},
  pages={265--290},
  year={1954},
  publisher={JSTOR}
}

@inproceedings{nedic2018improved,
  title={Improved convergence rates for distributed resource allocation},
  author={Nedi{\'c}, Angelia and Olshevsky, Alex and Shi, Wei},
  booktitle={2018 IEEE Conference on Decision and Control (CDC)},
  pages={172--177},
  year={2018},
  organization={IEEE}
}

@article{zhang2024geometric,
  title={Geometric Analysis of Matrix Sensing over Graphs},
  author={Zhang, Haixiang and Chen, Ying and Lavaei, Javad},
  journal={Advances in Neural Information Processing Systems},
  volume={36},
  year={2024}
}

@article{wang2016fully,
  title={A fully-decentralized consensus-based ADMM approach for DC-OPF with demand response},
  author={Wang, Yamin and Wu, Lei and Wang, Shouxiang},
  journal={IEEE Transactions on Smart Grid},
  volume={8},
  number={6},
  pages={2637--2647},
  year={2016},
  publisher={IEEE}
}

@article{van2014dc,
  title={DC power flow in unit commitment models},
  author={Van den Bergh, Kenneth and Delarue, Erik and D’haeseleer, William},
  journal={no. May},
  year={2014}
}

@article{yang2016state,
  title={A state-independent linear power flow model with accurate estimation of voltage magnitude},
  author={Yang, Jingwei and Zhang, Ning and Kang, Chongqing and Xia, Qing},
  journal={IEEE Transactions on Power Systems},
  volume={32},
  number={5},
  pages={3607--3617},
  year={2016},
  publisher={IEEE}
}

@article{kairouz2021advances,
  title={Advances and open problems in federated learning},
  author={Kairouz, Peter and McMahan, H Brendan and Avent, Brendan and Bellet, Aur{\'e}lien and Bennis, Mehdi and Bhagoji, Arjun Nitin and Bonawitz, Kallista and Charles, Zachary and Cormode, Graham and Cummings, Rachel and others},
  journal={Foundations and trends{\textregistered} in machine learning},
  volume={14},
  number={1--2},
  pages={1--210},
  year={2021},
  publisher={Now Publishers, Inc.}
}

@article{yang2019federated,
  title={Federated machine learning: Concept and applications},
  author={Yang, Qiang and Liu, Yang and Chen, Tianjian and Tong, Yongxin},
  journal={ACM Transactions on Intelligent Systems and Technology (TIST)},
  volume={10},
  number={2},
  pages={1--19},
  year={2019},
  publisher={ACM New York, NY, USA}
}

@article{nedic2010constrained,
  title={Constrained consensus and optimization in multi-agent networks},
  author={Nedic, Angelia and Ozdaglar, Asuman and Parrilo, Pablo A},
  journal={IEEE Transactions on Automatic Control},
  volume={55},
  number={4},
  pages={922--938},
  year={2010},
  publisher={IEEE}
}

@article{zhu2011distributed,
  title={On distributed convex optimization under inequality and equality constraints},
  author={Zhu, Minghui and Martinez, Sonia},
  journal={IEEE Transactions on Automatic Control},
  volume={57},
  number={1},
  pages={151--164},
  year={2011},
  publisher={IEEE}
}

@techreport{koloskova2024optimization,
  title={Optimization Algorithms for Decentralized, Distributed and Collaborative Machine Learning},
  author={Koloskova, Anastasiia},
  year={2024},
  institution={EPFL}
}

@article{beznosikov2023biased,
  title={On biased compression for distributed learning},
  author={Beznosikov, Aleksandr and Horv{\'a}th, Samuel and Richt{\'a}rik, Peter and Safaryan, Mher},
  journal={Journal of Machine Learning Research},
  volume={24},
  number={276},
  pages={1--50},
  year={2023}
}

@article{beznosikov2020derivative,
  title={Derivative-free method for composite optimization with applications to decentralized distributed optimization},
  author={Beznosikov, Aleksandr and Gorbunov, Eduard and Gasnikov, Alexander},
  journal={IFAC-PapersOnLine},
  volume={53},
  number={2},
  pages={4038--4043},
  year={2020},
  publisher={Elsevier}
}

@article{wang2022distributed,
  title={Distributed optimization with coupling constraints in multi-cluster networks based on dual proximal gradient method},
  author={Wang, Jianzheng and Hu, Guoqiang},
  journal={arXiv preprint arXiv:2203.00956},
  year={2022}
}

@article{liang2019distributed,
  title={Distributed smooth convex optimization with coupled constraints},
  author={Liang, Shu and Yin, George and others},
  journal={IEEE Transactions on Automatic Control},
  volume={65},
  number={1},
  pages={347--353},
  year={2019},
  publisher={IEEE}
}

@article{gong2023decentralized,
  title={Decentralized Proximal Method of Multipliers for Convex Optimization with Coupled Constraints},
  author={Gong, Kai and Zhang, Liwei},
  journal={arXiv preprint arXiv:2310.15596},
  year={2023}
}

@article{zhang2021distributed,
  title={Distributed convex optimization with coupling constraints over time-varying directed graphs},
  author={Zhang, Bingru and Gu, Chuanye and Li, Jueyou},
  journal={Journal of Industrial and Management Optimization},
  volume={17},
  number={4},
  pages={2119--2138},
  year={2021},
  publisher={Journal of Industrial and Management Optimization}
}

@article{li2018accelerated,
  title={Accelerated convergence algorithm for distributed constrained optimization under time-varying general directed graphs},
  author={Li, Huaqing and L{\"u}, Qingguo and Liao, Xiaofeng and Huang, Tingwen},
  journal={IEEE Transactions on Systems, Man, and Cybernetics: Systems},
  volume={50},
  number={7},
  pages={2612--2622},
  year={2018},
  publisher={IEEE}
}

@article{doan2017distributed,
  title={Distributed resource allocation on dynamic networks in quadratic time},
  author={Doan, Thinh T and Olshevsky, Alex},
  journal={Systems \& Control Letters},
  volume={99},
  pages={57--63},
  year={2017},
  publisher={Elsevier}
}

@article{salim2022dualize,
  title={Dualize, split, randomize: Toward fast nonsmooth optimization algorithms},
  author={Salim, Adil and Condat, Laurent and Mishchenko, Konstantin and Richt{\'a}rik, Peter},
  journal={Journal of Optimization Theory and Applications},
  volume={195},
  number={1},
  pages={102--130},
  year={2022},
  publisher={Springer}
}

@inproceedings{salim2022optimal,
  title={An optimal algorithm for strongly convex minimization under affine constraints},
  author={Salim, Adil and Condat, Laurent and Kovalev, Dmitry and Richt{\'a}rik, Peter},
  booktitle={International conference on artificial intelligence and statistics},
  pages={4482--4498},
  year={2022},
  organization={PMLR}
}

@article{gutknecht2002chebyshev,
  title={The Chebyshev iteration revisited},
  author={Gutknecht, Martin H and R{\"o}llin, Stefan},
  journal={Parallel Computing},
  volume={28},
  number={2},
  pages={263--283},
  year={2002},
  publisher={Elsevier}
}

@article{auzinger2011iterative,
  title={Iterative solution of large linear systems},
  author={Auzinger, Winfried and Melenk, J},
  journal={Lecture notes, TU Wien},
  year={2011}
}

@article{falsone2020tracking,
  title={Tracking-ADMM for distributed constraint-coupled optimization},
  author={Falsone, Alessandro and Notarnicola, Ivano and Notarstefano, Giuseppe and Prandini, Maria},
  journal={Automatica},
  volume={117},
  pages={108962},
  year={2020},
  publisher={Elsevier}
}

@article{wu2022distributed,
  title={Distributed optimization with coupling constraints},
  author={Wu, Xuyang and Wang, He and Lu, Jie},
  journal={IEEE Transactions on Automatic Control},
  volume={68},
  number={3},
  pages={1847--1854},
  year={2022},
  publisher={IEEE}
}

@article{chang2016proximal,
  title={A proximal dual consensus ADMM method for multi-agent constrained optimization},
  author={Chang, Tsung-Hui},
  journal={IEEE Transactions on Signal Processing},
  volume={64},
  number={14},
  pages={3719--3734},
  year={2016},
  publisher={IEEE}
}

@article{libsvm,
author = {Chang, Chih-Chung and Lin, Chih-Jen},
title = {LIBSVM: A library for support vector machines},
year = {2011},
issue_date = {April 2011},
publisher = {Association for Computing Machinery},
address = {New York, NY, USA},
volume = {2},
number = {3},
issn = {2157-6904},
url = {https://doi.org/10.1145/1961189.1961199},
doi = {10.1145/1961189.1961199},
journal = {ACM Trans. Intell. Syst. Technol.},
month = {may},
articleno = {27},
numpages = {27}
}

@inproceedings{dominguez2012decentralized,
  title={Decentralized optimal dispatch of distributed energy resources},
  author={Dominguez-Garcia, Alejandro D and Cady, Stanton T and Hadjicostis, Christoforos N},
  booktitle={2012 IEEE 51st IEEE conference on decision and control (CDC)},
  pages={3688--3693},
  year={2012},
  organization={IEEE}
}
